\documentclass[11pt,notitlepage,twoside]{article}
\usepackage{latexsym}
\usepackage{amsmath}
\usepackage{amsfonts}
\usepackage{amssymb}
\usepackage{amscd}
\renewcommand{\a }{\alpha }
\renewcommand{\b }{\beta }
\renewcommand{\d}{\delta }
\newcommand{\D }{\Delta }

\newcommand{\e }{\varepsilon }
\newcommand{\g }{\gamma}

\newcommand{\G }{\Gamma }
\renewcommand{\l }{\lambda }

\newcommand{\n }{\nabla }

\newcommand{\s }{\sigma }

\renewcommand{\o }{\omega }

\renewcommand{\O }{\Omega }

\newcommand{\ov}{\overline}
\newcommand{\intbar}{\mathop{\int\makebox(-13.5,0){\rule[4pt]{.7em}{0.3pt}}%
\kern-6pt}\nolimits}
\newcommand{\wtilde }{\widetilde}

\newcommand{\be}{\begin{equation}}
\newcommand{\ee}{\end{equation}}
\newcommand{\bes}{\begin{equation*}}
\newcommand{\ees}{\end{equation*}}
\newcommand{\ba}{\begin{eqnarray}}
\newcommand{\ea}{\end{eqnarray}}
\newcommand{\bas}{\begin{eqnarray*}}
\newcommand{\eas}{\end{eqnarray*}}
\newenvironment{pf}{\noindent{ \bf  Proof}.\enspace}{\rule{2mm}{2mm}\medskip}
\newenvironment{pfn}{\noindent{\bf  Proof}}{\rule{2mm}{2mm}\medskip}

\newcommand{\R}{\mathbb{R}}

\newcommand{\N}{\mathbb{N}}
\renewcommand{\o }{\omega }

\begin{document}

\author{Mohameden AHMEDOU$^{a}$\thanks{  E-mail: \texttt{Mohameden.Ahmedou@math.uni-giessen.de}}, Mohamed BEN AYED$^{b,c}$\thanks{E-mails: \texttt{M.BenAyed@qu.edu.sa and Mohamed.Benayed@fss.rnu.tn}} and
 Khalil EL MEHDI$^{b,d}$\thanks{ E-mail : \texttt{K.Jiyid@qu.edu.sa} } \\
 {\footnotesize
a : Department of Mathematics, Giessen University, Arndtstrasse 2, 35392, Giessen, Germany.}\qquad\quad\quad\\
{\footnotesize
b : Department of Mathematics, College of Science, Qassim University, Buraydah 51452, Saudi Arabia.}\\
{\footnotesize
c :  Facult{\'e} des Sciences de Sfax, Universit\'e de Sfax, Route Soukra 3000, Sfax, Tunisia.}\qquad\quad\qquad\qquad\\
{\footnotesize
 d : Facult\'e des Sciences et Techniques, Universit\'e de Nouakchott, Nouakchott, Mauritania.}\quad \qquad 
\quad}

\date{}

\title{\bf     Nirenberg problem on high dimensional  spheres: \\
Blow up with residual mass phenomenon \\
}

\newtheorem{lem}{Lemma}[section]
\newtheorem{pro}[lem]{Proposition}
\newtheorem{thm}[lem]{Theorem}
\newtheorem{rem}[lem]{Remark}
\newtheorem{cor}[lem]{Corollary}
\newtheorem{df}[lem]{Definition}

\maketitle

\noindent{\bf Abstract:} 
In this paper, we extend the analysis of the subcritical approximation of the Nirenberg problem on spheres recently conducted in \cite{MM19, MM}. Specifically, we delve into the scenario where the sequence of blowing up solutions exhibits a non-zero weak limit, which necessarily constitutes a solution of the Nirenberg problem itself. Our focus lies in providing a comprehensive description of such blowing up solutions, including precise determinations of blow-up points and blow-up rates. Additionally, we compute the topological contribution of these solutions to the difference  in topology between the level sets of the associated Euler-Lagrange functional. Such an  analysis is intricate due to the potential degeneracy of the involved solutions. We also provide a partial converse, wherein we construct blowing up solutions when the weak limit is  non-degenerate.

 \medskip
 


\section{Introduction and main Results}

Given a smooth positive function $K$ defined on the standard sphere $\mathbb{S}^n$, where $n \geq 3$, equipped with its standard metric $g_0$, the Nirenberg problem aims to determine a Riemannian metric $g$ conformally equivalent to $g_0$ such that the scalar curvature $R_g = K$. This can be formulated as solving the following nonlinear problem involving the critical Sobolev exponent:
$$
(\mathcal{N}_K) \qquad L_{g_0} u \, = \, K u^{(n+2)/(n-2)}, \quad u > 0 \quad\mbox{on}\quad \mathbb{S}^n,
$$
where  $ L_{g_{0}}:= -\D_{g_{0}}  \, + \, {n(n-2)}/{4}  $ denotes the conformal Laplacian.\\

This problem has garnered significant attention over the past half-century, with extensive research contributions.. See \cite{AH91, Bahri-Invariant, BC2, BCCH, CY, CGY, CL, CL1, CL2, Chen-Xu, Hebey, KW1, yyli1, yyli2, SZ} and the references therein.\\
Due to Kazdan-Warner topological obstructions  \cite{BEZ1, KW1}, the Nirenberg problem is not  solvable for every function $K$ . Thus, the central inquiry revolves around identifying sufficient conditions on $K$ for the problem to be solvable.\\
Regarding existence results, Euler-Poincaré type criteria have been established by Bahri-Coron \cite{BC2} and by Chang-Gursky-Yang \cite{CGY} on $\mathbb{S}^3$, and by Ben Ayed et al. \cite{BCCH} and Yanyan Li \cite{yyli2} on $\mathbb{S}^4$. Furthermore, these results have been extended to higher-dimensional spheres under various conditions, including flatness conditions near the critical point of $K$ \cite{yyli1}, in the perturbative framework \cite{CY}, or under pinching conditions \cite{Chen-Xu} and \cite{Malchiodi-Mayer}.

\noindent
 The Nirenberg problem has a variational structure. Indeed its solutions  correspond to the positive critical points of the functional
\be\label{functional}
I_{K}(u) := \frac{1}{2}  \| u \|^2 - \frac{ 1 }{ p+1}  \int_{ \mathbb{S}^n} K  | u |^{p+1} , \quad  \mbox{ with } p:= \frac{n+2}{n-2} ,
\ee
defined on  $H^1(\mathbb{S}^n)$  equipped with the scalar product
$$ 
\langle u,w\rangle := \int_{\mathbb{S}^n}  \n u \n w +\frac{n(n-2)}{4}
\int_{\mathbb{S}^n}  u w $$
and its associated norm
$$ \| u \|^2= \int_{\mathbb{S}^n} \mid \n u\mid^2 +\frac{n(n-2)}{4}
\int_{\mathbb{S}^n}  u^2.
$$
The functional $I_K$ fails to satisfy the Palais-Smale condition  and  the reason for such a lack of compactness is the existence of almost solutions of the equation $(\mathcal{N}_K)$. These almost solutions called \emph{bubbles} are defined as follows:
\begin{equation}\label{eq:deltatilde}
\wtilde{\d}_{(a,\l)}(x)=c_0
\frac{\l^{(n-2)/{2}}}{\bigl(2+(\l^2-1)(1-\cos d_{g_0}(x,a))\bigr)^{(n-2)/{2}}}, \quad \mbox{ with }\, \,  a\in \mathbb{S}^n \mbox{ and } \, \, \l>0,
\end{equation}
 where $d_{g_0}$ is the geodesic distance on $(\mathbb{S}^n,g_0)$ and $c_0=\left(n(n-2)\right)^{(n-2)/4}$.\\
 After performing a stereographic projection $\Pi$ with the
 point $-a$ as pole, the function $\wtilde{\d}_{(a,\l)}$ is
 transformed into

 \begin{equation}\label{eq:delta}
   \d_{(0,\l) } (y) = c_0\frac{\l^{(n-2)/{2}}}{(1+\l^2\mid
 y\mid^2)^{(n-2)/{2}}}, \qquad y \in \R^n,
 \end{equation}
 which is a solution of the problem \cite{CGS}
$$
 -\D u= u^\frac{n+2}{n-2} ,\, u > 0\,\quad \mbox{ in } \quad \R^n.
$$

One  way to overcome the lack of compactness of the functional $I_K$, which goes back to Yamabe \cite{Yamabe}, is to lower the exponent and first consider the subcritical approximation
 $$
(\mathcal{N}_{K,\tau}) \qquad L_{g_0} u \, = \, K u^{((n+2)/(n-2))-\tau}, \quad u > 0 \quad\mbox{on}\quad \mathbb{S}^n , 
$$
and its associated Euler-Lagrange functional
$$
 I_{K,\tau}(u) := \frac{1}{2}  \| u \|^2 - \frac{ 1 }{ p+1-\tau}  \int_{ \mathbb{S}^n} K  | u |^{p+1-\tau} , \quad u \in H^1(\mathbb{S}^n) \mbox{ with } p:= \frac{n+2}{n-2}.
 $$
Thanks to elliptic estimates either the solution $u_{\tau}$ remains uniformly bounded as   $\tau \to 0$ and hence converges strongly to a solution $ \o $ of the Nirenberg problem $(\mathcal{N}_K)$ or it blows up. In the latter  case,  following Schoen \cite{Schoen}, Yanyan Li \cite{yyli1, yyli2}, Chen-lin \cite{CL1,CL2} or Druet-Hebey-Robert \cite{DHR}, one performs a refined blow up analysis. 
It follows from such a blow up analysis that on $\mathbb{S}^3$ only single blow up can occur, whereas on $\mathbb{S}^4$ multiple blow up can occur, yet they are isolated simple, see \cite{BC2, CGY,yyli1,yyli2}.  However, the scenario shifts notably in higher-dimensional spheres. While on $\mathbb{S}^5, \mathbb{S}^6$ blowing up  solutions have finite energy \cite{CL2},  Chen-Lin \cite{CL} constructed  a sequence blowing up solutions of infinite energy. A fact that underscores the profound challenge presented by higher-dimensional spheres.  Recentyl Malchiodi-Mayer \cite{MM19, MM} have shown, that finite energy blowing up solutions having zero weak limit are all isolated simple. 
Motivated by   the intricate question of multiplicity of solutions to the Nirenberg problem, we undertake in this paper a systematic analysis of the scenario where a sequence of energy bounded solution of $(\mathcal{N}_{K,\tau})$ has a non zero weak limit. We point out that such a weak limit is necessarily a solution of $(\mathcal{N}_K)$ and 
it is known that such a  situation does not occur if the dimension of the sphere is less than six, see  \cite{BCCH, yyli1, yyli2, Bahri-Invariant}. For spheres of dimension $ n \geq 7$ our main result can be stated as follows:

\begin{thm}\label{th:blowup} Let $n \geq 7$, $0 < K \in C^3({\mathbb{S}^n})$ having only non-degenerate critical points $y_1, \cdots,y_q$ with $\D K(y_i) \neq 0$ for $1 \leq i \leq q$.
Let $(u_\tau)$ be a sequence of energy bounded solutions of $(\mathcal{N}_{K,\tau})$ converging weakly but not strongly $u_\tau \rightharpoonup \o \neq 0$. Then $\o$ is a solution of $(\mathcal{N}_K)$. Furthermore there exists $ N \in \N$ such that  $u_{\tau}$  has to blow up and takes the following form
$$
u_\tau \, = \,  \o \, + \,   \sum_{i = 1}^N   \frac{1}{K(a_{i, \tau})^{(n-2)/4}} \wtilde{\d}_{a_{i, \tau},\l_{i, \tau}} + v_\tau,
$$
where $\tilde{\delta}$  is defined in \eqref{eq:deltatilde} and $v_{\tau} \to 0$ in $H^1(\mathbb{S}^n)$.\\
Furthermore we have that:
\begin{enumerate}
  \item[i)]
  $a_{i, \tau}$  converges to a critical point  $y_i$ of $K$ with  $\D K(y_i) < 0$. 
  \item[ii)]
  $ \l_{i, \tau}\to \infty$, $\tau \ln\l_ {i, \tau} \to 0$ and  $\l_{i, \tau} d_{g_0} (a_{i, \tau},y_i) \to 0$ as $\tau \to 0.$
  \item[iii) ]
  For $i\ne j,$  \quad  $d_{g_0} ( a_{i, \tau},a_{j, \tau}) \, \geq \,   \frac{1}{2}\min\{ d_{g_0}(y_k, y_\ell ): y_k \neq y_\ell \in \mathcal{K}^\infty \}$ where \be\label{Kinfty} \mathcal{K}^\infty  := \{ y \in \mathbb{S}^n : \n K(y) =0 \mbox{ and } \D K(y) < 0\}. \ee
  \item[iv)]
 There exists  a dimensional constant $\kappa_1(n) > 0$ such that
  \be\label{lll}
   \frac{1}{\l_{i, \tau} ^2} = - \kappa_1(n) \frac{K(y_i) }{ \D K(y_i)} \, \tau (1+o_{\tau}(1)).
   \ee
\end{enumerate}
In addition, if $\o$ is a non-degenerate critical point of  $I_K$, then  $u_{\tau}$ is a non-degenerate critical point of  $I_{K,\tau}$ and
  $$ \mbox{Morse index} (I_{K,\tau}, u_{\tau}) =  (N+1) + \mbox{Morse index} (I_{K} , \o ) +  \sum_{i=1}^N ( n - \mbox{Morse index} ( K, y_i )).$$
\end{thm}

Regarding the proof strategy of Theorem \ref{th:blowup}, it is essential to provide some insights. While conventional blow-up analysis techniques typically rely heavily on precise pointwise $C^0$-estimates of the blowing-up solutions $u_{\tau}$ and extensively employ Pohozaev identities \cite{CL1, CL2, DHR, KMS09, yyli1, yyli2, MM19, Schoen}, our approach diverges from this path.
Our strategy hinges on deriving balancing conditions governing the parameters of concentration, utilizing refined asymptotic estimates of the gradient within the so-called "neighborhood at infinity". These balancing conditions emerge through testing the equation with vector fields representing the leading terms of the gradient concerning the concentration parameters. By analyzing these conditions, we can extract all pertinent information regarding the blow-up phenomenon. 
We highlight that the potential degeneracy of the weak limit introduces an additional challenging dimension to the problem. Estimating parameters within the nontrivial kernel of the linearized operator becomes particularly arduous in such scenarios.
It is worth noting that circumventing the reliance on pointwise estimates and Pohozaev identities can be advantageous, especially in the exploration of non-compact variational problems where non-simple blow-ups may arise, such as in the singular mean-field equation with quantized singularities \cite{BT, KLin, WZ, WZb, DW}, and in the Nirenberg problem on half-spheres \cite{AB20b}. The presence of non-simple blow-up points significantly complicates the task of establishing pointwise $C^0$-estimates, making our method particularly valuable in such contexts.

\bigskip
\noindent
In the next theorem we provide, under generic condition, the following  converse of  Theorem \ref{th:blowup}.
\begin{thm}\label{th:blowupConv}
Let $n \geq 7$, $0 < K \in C^3({\mathbb{S}^n})$. Let $y_1, \cdots, y_N$ be distinct non-degenerate critical points of $K$ with $\D K(y_i) <0$ and $\o$ be a non-degenerate solution of $(\mathcal{N}_K)$. Then, as $\tau\to 0$, there exists a  unique solution $u_{\tau,\o, y_1,...,y_N}$ of $(\mathcal{N}_{K,\tau})$  satisfying: $u_{\tau,\o, y_1,...,y_N} -\o$ develops exactly one bubble at each point $y_i$ and converges weakly to zero in $H^1(\mathbb{S}^n)$ as $\tau \to 0$. More precisely there exist $\l_{1,\tau}$,..., $\l_{N,\tau}$ having the same order as $\tau^{-1/2}$ and points $a_{i,\tau}\to y_i$ for all $i$ such that
\begin{align*}
& \bigg |\bigg | u_{\tau,\o, y_1,...,y_N} -\o-\sum_{i=1}^NK(a_i)^{\frac{2-n}{4}}\wtilde{\d}_{a_{i,\tau}, \l_{i,\tau}} \bigg |\bigg | \to 0\quad \mbox{and}\\
&  I_\tau (u_{\tau,\o, y_1,...,y_N}) \to \frac{1}{n}  \vert \vert \o\vert\vert^{ 2 } + \frac{1}{n}  S_n \left(\sum_{i=1}^N K(a_i)^{\frac{2-n}{2}}\right),
\end{align*}
as $\tau\to 0$,  where
\begin{equation}\label{eq:sn} S_n := c_0 ^{2n/(n-2)} \int_{\R^n} \frac{1}{(1+ | x | ^2 )^n } dx.\end{equation}
In addition,  $u_{\tau,\o, y_1,...,y_N}$ is non-degenerate for $I_{K,\tau}$ and
 $$ \mbox{Morse index} (I_{K,\tau}, u_{\tau,\o, y_1,...,y_N}) =  (N+1) + \mbox{Morse index} (I_{K} , \o ) +  \sum_{i=1}^N ( n - \mbox{Morse index} ( K, y_i )).$$
\end{thm}
We would like to underscore that our main findings in this paper, namely Theorems \ref{th:blowup} and \ref{th:blowupConv}, serve as fundamental components in our subsequent work \cite{ABE}. Specifically, they play a pivotal role in demonstrating the existence of infinitely many solutions to Nirenberg's problem under perturbative conditions.

\bigskip

The subsequent sections of the paper are structured as follows:
In Section 2, we perform  a finite-dimensional reduction while 
Section 3 delves into providing a meticulous expansion of the gradient within the neighborhood at infinity.
The refined blow-up analysis of finite energy approximated solutions to problem $(\mathcal{N}_{K,\tau})$, with a non-zero weak limit, is detailed in Section 4 and Section 5 is dedicated to the proof of Theorem \ref{th:blowupConv}.
Lastly, in the appendix, we gather some indispensable estimates crucial for substantiating various assertions throughout this paper.

\section{ Finite dimensional reduction}

  Let $\o$ be a solution of $(\mathcal{N}_K)$  and let  $N_0(\o)$ be the kernel of the associated quadratic form
defined by:
\begin{equation}\label{qw} Q_\o (h):= \| h \| ^2 - \frac{n+2}{n-2} \int_{\mathbb{S}^n} K \o ^{4/(n-2)} h^2 \qquad \mbox{ for } h \in H^1(\mathbb{S}^n). \end{equation}
 Let $m_2$ be the dimension of $N_0(\o)$ and  $(e_1,\cdots,e_{m_2})$  be an orthonormal basis of $N_0(\o)$.  We set
  $$
H_0(\o):= \, span(\o) \oplus span(e_1, \cdots,e_{m_2}).
$$
Following  M. Mayer,  see Lemma 3.6 and Proposition 3.7 in \cite{Mayer-Thesis}, we parameterize  a neighborhood of $\o$ by
\be \label{ualphbeta} u_{\a,\b} := \a (\o + \sum_{i=1}^{m_2}  \b_i e_i + h(\b))  \quad  \mbox{ with }  h(\b) \perp H_0(\o); \, h(\beta) = O(||\b||^2) \mbox{ and } \| h(\b) \|_{C^2} \to  0, \ee
where $\b:= (\b_1, \cdots , \b_{m_2}) \in \R^{m_2}$, $\a$ close to $1$ and  the function  $u_{\a,\b} $ satisfies
\be \label{ualphbeta2} \langle  \n I_K(u_{\a,\b}), h \rangle := \langle u_{\a,\b},h\rangle  - \int _{ \mathbb{S}^n} K u_{\a,\b}^{p} h = 0  \quad \mbox{ for each } h\in H_0(\o)^{\perp} .\ee

Next for $\o$ a solution of $(\mathcal{N}_K)$ whose Kernel is of dimension $m_2$,  $N \in \N_0$ and $\mu$ a small positive constant, we define the so called \emph{neighborhood at infinity} $V(\o, N,\mu)$ as follows:

\begin{align}\label{neigh-at-Infty}
{V}(\o, N,\mu) :=  \Big \{  & u \in H^1(\mathbb{S}^n):  \,
 \exists \,  \l_1, \cdots, \l_N > {\mu^{-1}} \, ; \,   \,  \exists \, a_1, \cdots,a_{N} \in { \mathbb{S}^n}, \quad  \mbox{ with } \nonumber \\
 &  \e_{ij} < \mu ; \, \,   \exists \, \a_0 \in(1-\mu,1+\mu) ;\, \,  \exists \, \,  \b \in \R^{m_2} \, \, \,   \mbox{ with }\, \,  \| \b \| \leq c\, \mu  \,  \nonumber \\
   &   \mbox{ such that } \| u -  \sum_{i = 1}^N K(a_i)^{(2-n)/4} \wtilde{\d}_{a_i,\l_i} - u_{\a_0,\b} \| < \mu  \Big \},
\end{align}
where
\begin{equation}\label{eppe}
\e_{ij}:= \Big( \frac{\l_i}{\l_j} + \frac{\l_j}{\l_i} + \frac{1}{2} \l_i \l_j (1 - \cos  d (a_i,a_j))\Big)^{(2-n)/2}.
\end{equation}

Following Bahri-Coron \cite{BCd}  we consider, for  $u \in V(\o, N ,\mu)$, the following minimization problem

\begin{equation}\label{eq:min}
  Min  \left   \{  \|  u  \, - \, \sum_{i=1}^{N}\a_i  \wtilde{\d}_{a_i,\l_i} - u_{\a_0,\b} \|;  \a_i > 0 ; \b\in \R^{m_2} ;  \l_i > 0,  \, a_i \in \mathbb{S}^n  \right  \}.
\end{equation}
We then have the following proposition whose proof is  identical, up to minor modification, to  the one  of  Proposition 7 in \cite{BCd}.

\begin{pro}\label{p:min}
For any $N  \in \N_0$ there exists $\mu_0 > 0$ such that if $\mu < \mu_0$ and $u \in V( \o, N , \mu)$ the minimization problem \eqref{eq:min}  has, up to permutation of the indices, a unique solution.
\end{pro}
Hence it follows from Proposition \ref{p:min} that every $u \in V(\o,  N , \mu)$ can be written in a unique way as
\begin{align}
& u \, = \, \sum_{i=1}^{N} \a_i  \wtilde{\d}_{a_i,\l_i} \, + \, u_{\a_{0},\b} \, +  \, v, \qquad \mbox{ where }\label{F1}\\
& a_i \in  \mathbb{S}^n , \, i \leq N , \quad  \a_i^{4/(n-2)} K(a_i) = 1 + o(1) \, \,  \forall \, i\geq 1,  \quad \a_0 = 1 + o(1),  \label{alphai0}
\end{align}
and  $v \in H^1(\mathbb{S}^n)$ satisfying
\begin{align}
&   \| v\| <  \mu, \quad <v, \psi> = 0, \mbox{ for } \psi \in  E_{\o, a,\l} ^\perp \mbox{ where } \label{eq:V0} \\
& E_{\o, a,\l}:=\mbox{span} \{   \wtilde{\d}_i, \frac{\partial  \wtilde{\d}_i}{\partial \l_i}, \frac{\partial  \wtilde{\d}_i}{\partial a_i} , \, u_{\a_0,\b} , \frac{\partial u_{\a_0,\b}}{\partial \b_k} ;   1 \leq i\leq N; \,  ; \, k\leq m_2\} , \label{Ewal}
\end{align}
where $  \wtilde{\d}_i :=  \wtilde{\d}_{a_i,\l_i}$ .
Furthermore a  refined Struwe energy type decomposition (\cite{Struwe}) has been proved in \cite{MM19}. Namely the following result holds: 

\begin{pro}\label{lambdaepsilonw}[\cite{MM19}, Proposition 3.1]
Let $ u_\tau$  be an energy bounded solution of $(\mathcal{N}_{K,\tau})$ which blows up. We assume that  there exists a positive solution $\o$ of $(\mathcal{N}_K)$ such that $u_\tau \rightharpoonup \o$ (but $u_\tau \nrightarrow \o$). Then there exists $ N $ such that $u_\tau$ can be written as

 \begin{equation}\label{kkk111} u_\tau:= u_{\a,\b} + \sum_{ i = 1}^N \a_i  \wtilde{\d}_{ a_i, \l_i} + v_\tau \in V(\o, N ,\mu )\quad \mbox{ with } v_\tau \in E_{\o,a,\l}^\perp \mbox{ and } \tau \ln \l_i = o_\tau (1) \, \, \forall \, i. \end{equation}
\end{pro}
In what follows, we assume that $u_\tau$ is decomposed as in \eqref{kkk111}.
Now, we are going to estimate the remainder term $v_\tau $  in Proposition \ref{lambdaepsilonw}.
To this aim, we need to give some information about the second variation of the Euler-Lagrange functional with respect to the infinite dimensional variable $v$. First, for  a solution $\o$  of $(\mathcal{N}_K)$ (not necessary  a non-degenerate one),  we decompose $H^1(\mathbb{S}^n)$ as follows:
\be \label{vvv*1}  H^1(\mathbb{S}^n) := N_-(\o) \oplus H_0(\o) \oplus N_+(\o) \, \, \quad \mbox{where} \, \, H_0(\o):= span\{\o\} \oplus N_0(\o)\ee
where  $N_-(\o)$, $N_0(\o)$ and $ N_+(\o)$ are respectively the space of negativity, of nullity and of positivity  of the quadratic form $Q_\o$ (defined by \eqref{qw}) in $span\{\o\} ^\perp$. Note that these spaces are orthogonal spaces with respect to $\langle.,.\rangle$ and the bilinear form $B_\o(.,.)$ ($:= \int_{\mathbb{S}^n}  K \o^{p-1} . . $). Furthermore, the sequence of the eigenvalues (denoted by $(\s_i)$) corresponding to $Q_\o$ satisfies $ \s_i \nearrow 1$. Therefore, there exists a constant $c >0$ such that
\begin{equation}\label{qw+-}
Q_\o(h) \leq -c \, \| h \| ^2 \quad \mbox{ for each } h \in N_-(\o) \quad ; \quad Q_\o(h) \geq c \, \| h \| ^2 \quad \mbox{ for each } h \in N_+(\o).\end{equation}

\begin{lem}\label{Qwalpositive} Assume that the $\e_{ij}$'s are small. Let $v \in E_{\o,a,\l} ^\perp$ (defined in \eqref{Ewal}).  Written $ v $ as 
$$ v := v_- + v_0 + v_+ \quad  \mbox{ with } \quad
 v_- \in N_-(\o) \, \, ; \, \,  v_0 \in H_0(\o) \, \, ; \, \, v_+ \in N_+(\o) \quad (\mbox{see } \eqref{vvv*1}). $$  
 Let us define 
 \be \label{ff7} Q_{\o,a,\l} (v) := \| v \|^2 - p \sum \int_{\mathbb{S}^n} \wtilde{\d}_i ^{p-1} v ^2 - p  \int_{\mathbb{S}^n} K \o ^{p-1} v ^2. \ee
 Then, $ Q_{\o,a,\l} $ is a non-degenerate quadratic form on $ E_{\o, a , \l }^\perp $ . More precisely, there exists a positive constant $\underline{c} $ such that
\begin{align}  & \| v_0 \| = o( \| v \|  ),  \label{0q0} \\
& Q_{\o,a,\l} (v_- ) = Q_\o (v_-) + o( \| v_- \| ^2 ) \leq -\,  \underline{c}\,  \| v_-\|^2, \label{0q-}\\
& Q_{\o,a,\l} (v_+ ) \geq  \underline{c} \,  \| v_+\|^2 + o( \| v \| ^2), \label{0q+}
\end{align}
where $Q_\o$ is defined in  \eqref{qw}.
\end{lem}

\begin{pf} First, we remark that, \eqref{0q0}-\eqref{0q+} imply that $ Q_{\o,a,\l} $ is non-degenerate on $ E_{\o,a,\l} ^\perp $. \\
Second,  note that  the spaces  $H_0(\o)$, $N_-(\o)$ and $N_+(\o)$  are orthogonal spaces with respect to $\langle.,.\rangle$ and the associated bilinear form $B_\o(.,.)$ ($:= \int_{\mathbb{S}^n}  K \o^{p-1} . . $).

We start by proving \eqref{0q0}. Since $v_0 \in H_0(\o)$, it follows that $ v_0 = \g_0 \o + \sum \g_i e_i$ where $(e_1,\cdots,e_{m_2})$ is an orthonormal basis of $N_0(\o)$. Using the fact that $v \in E_{\o,a,\l}^\perp$ (which implies that $ v \perp u_{\a,\b}$ and $v \perp \partial u_{\a,\b} / \partial \b_i$ for each $i $), it follows that
\begin{align*} & \g_0 =  \langle v_0 , \o \rangle  =  \langle v , \o \rangle = (1/ \a) \langle v , u_{\a,\b}  \rangle  -  (1/\a) \sum \b_i \langle v , e_i \rangle -  \langle v , h(\beta) \rangle = o( \| v \| ), \\
& \g_i = \langle v_0 , e_i \rangle = \langle v , e_i \rangle = \langle v , e_i + (\partial h(\b)/ \partial \b_i)  \rangle - \langle v ,  \partial h(\b)/ \partial \b_i  \rangle = o( \| v \| ) \quad \forall \,\,   1 \leq i \leq m_2 \end{align*}
(by using the smallness of $\b$ and $h(\b)$ in the $C^1$ sense with respect to $\b$). This ends the proof of \eqref{0q0}.

Concerning \eqref{0q-}, we have
$$ Q_{\o,a,\l} (v_- ) = \| v_- \| ^2 - p \sum \int \wtilde{ \d}_i ^{p-1} v_- ^2 - p \int K \o ^{p-1} v_- ^2 = Q_\o ( v_-) - p \sum \int  \wtilde{ \d}_i ^{p-1} v_- ^2 .$$
Observe that, since $v_-$ belongs to a  fixed finite dimensional space, we derive that $\| v_- \|_{\infty} \leq c \| v_- \| $ and therefore
$$ \int \wtilde{ \d}_i ^{p-1} v_-^2 \leq \| v_- \|^2_{\infty}  \int \wtilde{\d}_i ^{p-1}  = o(  \| v_- \|^2 ) \quad \mbox{ for each } i.$$
Hence the proof of \eqref{0q-} follows by using \eqref{qw+-}. \\
It remains to prove \eqref{0q+}.
Note that, using Proposition 3.1 of \cite{B1}, there exists a constant $\underline{c}_1> 0$ such that
\begin{equation}\label{h1}
\| h \| ^2 - \frac{n+2}{n-2}  \sum_{i=1}^N \int _{\mathbb{S}^n} \wtilde{\d}_{a_i, \l_i} ^{\frac{4}{n-2}} h^2  \geq \underline{c}_1 \| h \| ^2 \quad \mbox{ for each } h \in E_{a,\l}^\perp.
\end{equation} where $E_{a,\l}^\perp$ is introduced in  \eqref{Ewal} (with $\o=0$). \\
In addition, the sequence of the eigenvalues (denoted by $(\s_i)$)  corresponding to $Q_\o$ (defined by \eqref{qw}) satisfies $\s_i \nearrow 1$. Let $N_k(\o)$ be the eigenspace associated to the eigenvalue $\s_k$. These spaces are orthogonal with respect to $\langle ., .\rangle$ and the bilinear form $B_\o$. Let $\s_{k_0}:= \min\{ \s_i : \s_i > 0\}$. Hence it is easy to see that $N_+(\o)= \oplus_{ k \geq k_0} N_k(\o)$. Furthermore, for $ k \in \N$, it holds
\begin{equation}\label{h2} Q_\o (h) \geq \s_k \| h \| ^2 \qquad \mbox{ for each } h \in \oplus_{ j \geq k } N_j(\o). \end{equation}
Let $k_1$ be such that $\s_{k_1} \geq 1- \underline{c}_1 /4 $. We decompose $N_+(\o)$ as follows:
$$ N_+(\o) := (\oplus_{k_0 \leq k \leq k_1} N_k(\o)) \oplus (\oplus_{k > k_1} N_k(\o)) := N_+^{0,1}(\o) \oplus N_+^{1}(\o).$$
Note that $ N_+^{0,1}(\o)$ is a fixed finite dimensional  space.
Now, since  $v_+ \in  N_+(\o)$, it holds that
\be\label{qaz2} v_+:= v_0^+ + v_1^+ \qquad \mbox{ where } v_0^+ \in  N_+^{0,1}(\o) \mbox{ and } v_1^+ \in  N_+^1(\o).\ee
Hence, using the orthogonality of the spaces $N_+^{0,1}(\o)$ and $ N_+^1(\o)$, it follows that
\begin{align*} Q_{\o,a,\l} (v_+) = \| v_0^+ \| ^2 +  \| v_1^+ \| ^2 & - p  \sum_{i=1}^N \int _{\mathbb{S}^n} \d_{a_i, \l_i} ^{\frac{4}{n-2}} \{ (v_0^+)^2 + (v_1^+)^2 + 2 v_0^+ v_1^+ \} \\
 & - p  \int _{\mathbb{S}^n} K \o ^{\frac{4}{n-2}} \{ (v_0^+) ^2 + (v_1^+)^2 \}. \end{align*}
Observe that
\begin{align*}
& \| v_0^+ \| ^2 - \frac{n+2}{n-2}  \int _{\mathbb{S}^n} K \o ^{\frac{4}{n-2}} (v_0^+) ^2 = Q_\o(v_0^+)  \geq \s_{k_0} \| v _0^+ \| ^2  \quad \mbox{ (by using \eqref{h2}) } \\
&  \int _{\mathbb{S}^n} \wtilde{\d}_{a_i, \l_i} ^{\frac{4}{n-2}} (v_0^+)^2 \leq \| v_0^+\|^2_\infty \int _{\mathbb{S}^n} \wtilde{\d}_{a_i, \l_i} ^{\frac{4}{n-2}}  = o( \| v_0^+ \| ^2)\\
&  \int _{\mathbb{S}^n} \wtilde{\d}_{a_i, \l_i} ^{\frac{4}{n-2}} | v_0^+ | | v_1^+ | \leq \| v_0^+\|_\infty \int _{\mathbb{S}^n} \wtilde{\d}_{a_i, \l_i} ^{\frac{4}{n-2}}  | v_1^+ | =  o( \| v_0^+ \| \| v_1^+ \|) = o( \| v_1^+  \| ^2 + \| v_0^+ \|^2 )\\
& - p  \int _{\mathbb{S}^n} K \o ^{\frac{4}{n-2}}  (v_1^+)^2 = Q_\o (v_1^+) - \| v_1^+ \|^2 \geq (\s_{k_1} - 1 ) \| v_1^+ \| ^2 \geq -(\underline{c}_1/4) \| v_1^+ \| ^2
\end{align*}
where, for the last formula,  we have used  \eqref{h2} and the choose of $k_1$. Combining these estimates, we get
\be\label{qaz1} Q _{\o,a,\l} (v_+)  \geq \{ \| v_1^+ \| ^2  - p  \sum_{i=1}^N \int _{\mathbb{S}^n} \wtilde{\d}_{a_i, \l_i} ^{\frac{4}{n-2}} (v_1^+)^2 \} + \s_{k_0} \| v _0^+ \| ^2  -\frac{\underline{c}_1}{4} \| v_1^+ \| ^2 + o( \| v_1^+ \| ^2 + \| v_0^+ \| ^2) .\ee
Note that the function $v_1^+$ is not necessarily in $E_{a,\l}^\perp$. For this reason we write $ v_1^+ := \sum t_i \psi _i + \tilde{v}_1^+$ with  $\tilde{v}_1^+ \in E_{a,\l}^\perp$ where the $\psi_i$'s are the functions $ \wtilde{\d}_j$'s and their derivatives with respect to $ \l_j$ and $a_j^k$. Let $\psi _i \in \{ \wtilde{\d}_i,\l_i  \partial \wtilde{\d}_i / \partial \l_i  , (1/\l_i)  \partial \wtilde{\d}_i / \partial \a_i ^k  \}$, it follows that
\begin{align*} t_i + o(\sum | t_k |) & = \langle v_1^+ , \psi_i \rangle =  \langle v  , \psi_i \rangle -  \langle v_0 , \psi_i \rangle  - \langle v_- , \psi_i \rangle - \langle v_0^+ , \psi_i \rangle\\
&  = O \Big( \int \wtilde{\d}_i ^p ( | v_0| + | v_- | + | v_0^+ | ) \Big) = o(  \| v_0 \| + \| v_- \| + \| v_0^+ \|) = o(  \| v \|)\end{align*} 
by using the fact that these functions are in  fixed finite dimensional spaces.  Thus we derive that $ \| v_1^+ \|^2 = \| \tilde{v}_1^+ \|^2 +o( \| v \|^2 )$ and therefore, using \eqref{h1}, we get
\begin{align}  \| v_1^+ \| ^2  - p  \sum_{i=1}^N \int _{\mathbb{S}^n} \wtilde{\d}_{a_i, \l_i} ^{\frac{4}{n-2}} (v_1^+)^2  & = \| \tilde{v}_1^+ \| ^2  - p \sum_{i=1}^N \int _{\mathbb{S}^n} \wtilde{\d}_{a_i, \l_i} ^{\frac{4}{n-2}} (\tilde{v}_1^+)^2  + o( \| \tilde{v}_1^+ \| ^2  +\| v \| ^2 )\nonumber  \\
& \geq \frac{1}{2} \underline{c}_1  \| v_1^+ \| ^2 + o( \| v \| ^2 )   . \label{qaz3} \end{align}
Combining \eqref{qaz2}, \eqref{qaz1} and \eqref{qaz3}, we get
$$ Q _{\o,a,\l} (v_+)  \geq \frac{1}{4} \underline{c}_1   \|  v_1^+ \| ^2 +  \s_{k_0} \| v _0^+ \| ^2 + o( \| v_1^+ \| ^2 + \| v_0^+ \| ^2) + o( \| v \| ^2 )  \geq c \| v_+ \| ^2 + o( \| v \| ^2 ) .$$
Thus the result follows.
\end{pf}

Next, we prove the following proposition which gives the estimate of the infinite dimensional part of $ u_\tau $. 

\begin{pro}\label{eovvw}
Let $v_\tau $ be the  remainder term  defined in Proposition \ref{lambdaepsilonw}. Then there holds:
$$ \| v_\tau  \| \leq c R(\tau ,a,\l) \quad \mbox{ with }
R(\tau ,a,\l):= \tau  +  \sum_{i=1}^N \frac{| \n K(a_i) | }{\l _i} + \frac{1}{\l _i^2} +
 \sum \e _{ij }^{\frac {n+2}{2(n-2)}}(\ln \e _{ij}^{-1})^{\frac{n+2}{2n}} . $$
\end{pro}

 \begin{pf}
Using \eqref{vvv*1},  $v_\tau $ can be decomposed as  follows
\be\label{emna13} v_\tau := v_\tau^- + v_\tau^0 + v_\tau^+ \quad \mbox{ where } v_\tau ^0 \in  H_0(\o) \, \, ; \, \, v_\tau ^- \in  N_-(\o) \mbox{ and }  v_\tau ^+ \in  N_+(\o) .\ee
Since $v_\tau \in E_{w,a,\l}^\perp$, using Lemma \ref{Qwalpositive}, we get that $ v_\tau^0 = o( \| v_\tau \| )$. For the other parts, note that $v_\tau^+$ and $v_\tau^-$ are not necessarily in $E_{a,\l}^\perp$ but they are in $H_0(\o)^\perp$. \\
Now we will focus on estimating $v_\tau^-$. Since $ u_\tau $ is a solution of $ (\mathcal{N}_{K,\tau} ) $, we get 
$$ \langle \n I_{K,\tau}(u_\tau) , {v }_\tau ^- \rangle = 0 $$ which is equivalent to (by using  \eqref{nablaI} and  \eqref{ti6})
\begin{align}
 \sum \a_j \langle \wtilde{\d}_j , v_\tau^- \rangle +  \langle u_{\a,\b} , v_\tau^- \rangle + \| v_\tau^- \|^2 & = \int  K u_\tau ^{p-\tau} v_\tau^- \nonumber\\
& =   \int K \ov{u}_\tau ^{p-\tau} v_\tau^- + p \int K \ov{u} _\tau ^{p-1-\tau} v_\tau v_\tau^- + o(\| v_\tau\| \| v_\tau^-\|) \label{emna1}\end{align}
where $\ov{u}_\tau := u_\tau - v_\tau.$
Using \eqref{ti5}, \eqref{alphai0} and the fact that $v_\tau^-$ is in a finite dimensional space  (which implies that $\| . \|_\infty \leq c \| . \|$), we derive that
\begin{align}
 \int K \ov{u} _\tau ^{p-1-\tau} v_\tau v_\tau^- & = \sum_{i=1}^N \int K (\a_i \wtilde{\d}_ i) ^{p-1-\tau} v_\tau  v_\tau ^- + \int K u_{\a,\b}^{p-1-\tau} v_\tau  v_\tau ^- + o(\| v_\tau\| \| v_\tau ^-\|)\nonumber \\
 & = \sum_{i=1}^N O\big( | v_\tau ^- |_\infty \int  \wtilde{\d}_ i ^{p-1} | v_\tau | \Big)  +  \int K \o^{p-1} ( v_\tau ^0 +v_\tau ^- + v_\tau^+ ) v_\tau^- + o(\| v_\tau \| \| v_\tau ^-\|) \nonumber \\
& =  \int K \o^{p-1} ( v_\tau ^-  )^2 + o(\| v_\tau \| \| v_\tau ^-\|)  \label{emna2}
 \end{align} where we have used the orthogonality of $v_\tau ^-$, $v_\tau^0$ and $ v_\tau ^+$ with respect to $\int K\o^{p-1} .. $ . Thus \eqref{emna1} and  \eqref{emna2} imply that
\be\label{emna3}  - Q_\o (v_\tau^-) + o(\| v_\tau\| \| v_\tau ^-\|) =  \sum \a_j \langle \wtilde{\d}_j , v_\tau^- \rangle +  \langle u_{\a,\b} , v_\tau^- \rangle  - \int K \ov{u}_\tau ^{p-\tau} v_\tau^- := \ell (v_\e^-) . \ee
Observe that, using $ \| v_\tau^- \| _\infty \leq c \| v_\tau ^- \|$, it holds that
\be\label{emna11}
| \langle  \wtilde{\d}_i , v_\tau^- \rangle | \leq c  \int \wtilde{\d}_i ^p | v_\tau^- |  \leq c   \| v_\tau^- \|_\infty  \int \wtilde{\d}_i ^p  \leq c  \frac{ \| v_\tau^- \| }{\l_i ^{(n-2)/2}} \quad ; \quad  \int \wtilde{\d}_i  | v_\tau^- |  \leq c\frac{ \| v_\tau ^- \| }{\l_i ^{(n-2)/2}} . \ee
\begin{align} \int K \ov{u}_\tau ^{p-\tau} v_\tau^- & = \int K {u}_{\a,\b} ^{p-\tau} v_\tau^- + \sum O\Big( \int u_{\a,\b} ^{p-1} \d_i | v_\tau^-| + \int \d_i ^{p-\tau} | v_\tau^-| \Big) \nonumber \\
& = \int K {u}_{\a,\b} ^{p} v_\tau^-  +  O\Big( \| v_\tau^- \| \Big( \tau + \sum \frac{1}{\l_i ^{(n-2)/2}}\Big)\Big) . \label{emna4}\end{align}
Therefore, combining \eqref{qw+-}, \eqref{emna3}-\eqref{emna4} and \eqref{ualphbeta2}, we get
\be\label{emna5}  c \, \| v_\tau ^- \| ^2 \leq - Q_\o (v_\tau ^- ) \leq c \, \| v_\tau ^-  \| \Big( \tau + \sum \frac{1}{\l_i ^{(n-2)/2}} + o( \| v_\tau \| ) \Big) .\ee
It remains to estimate the $v_\tau ^+$-part. Recall that $ v_\tau ^+ \in N_+(\o)$ and it is not necessarily in $E_{a,\l}^\perp$. Note that Eq \eqref{emna1} holds with $v_\tau ^+$ instead of $v_\tau ^-$. Observe that, using \eqref{ti5} and \eqref{alphai0}, it holds
\begin{align}
 \int K \ov{u} _\tau ^{p-1-\tau} v_\tau v_\tau^+ & = \sum_{i=1}^N \int K (\a_i \wtilde{\d}_ i) ^{p-1-\tau} v_\tau  v_\tau ^+ + \int K u_{\a,\b}^{p-1-\tau} v_\tau  v_\tau ^ + + o(\| v_\tau\| \| v_\tau ^+\|)\nonumber \\
 & = \sum_{i=1}^N   \int  \wtilde{\d}_ i ^{p-1}  v_\tau   v_\tau ^+  +  \int K \o^{p-1} ( v_\tau ^0 +v_\tau ^- + v_\tau^+ ) v_\tau^+ + o(\| v_\tau\| \| v_\tau ^+\|) \nonumber \\
& =   \sum_{i=1}^N   \int  \wtilde{\d}_ i ^{p-1}   (v_\tau ^+)^2  + \int K \o^{p-1} ( v_\tau ^+  )^2 + o(\| v_\tau\| \| v_\tau ^+\|) , \label{emna22}
 \end{align}
 by using the fact that $ \| v_\tau^- \| _\infty \leq c \| v_\tau ^- \| \leq c \| v_\tau  \|$, $ \| v_\tau^0 \| _\infty \leq c \| v_\tau ^0 \| \leq c \| v_\tau  \| $ and the fact that  $\int K\o^{p-1} v_\tau ^-v_\tau ^+ = \int K\o^{p-1} v_\tau ^0v_\tau ^+=0 $ . Hence we obtain
 \be\label{emna9} Q_{\o,a,\l} (v_\tau^+) +  o(\| v_\tau\| \| v_\tau ^+\|) =  - \sum \a_j \langle \wtilde{\d}_j , v_\tau^+ \rangle -  \langle u_{\a,\b} , v_\tau^+ \rangle  +  \int K \ov{u} _\tau ^{p-\tau} v_\tau^+ := - \ell (v_\tau^+).\ee
Using Lemma \ref{Qwalpositive}, we derive that $Q_{\o,a,\l} (v_\tau^+) \geq c \| v_\tau^+\|^2 + o(\| v_\tau\|^2)$, hence it remains to estimate the linear part $\ell(v_\tau^+)$. In fact, using  \eqref{ti5}, we have
\begin{align} \int K \ov{u}_\tau ^{p-\tau} v_\tau^+ = & \sum \int_{\mathbb{S}^n} K ( \a_i \wtilde{\d}_i)^{p-\tau} v_\tau^+ + \int_{\mathbb{S}^n} K u_{\a,\b}^{p-\tau} v_\tau ^+ \nonumber\\
&  + O\Big(  \sum_{i\neq j} \int  (\wtilde{\d}_i \wtilde{\d}_j)^{p/2} | v_\tau ^+|  + \sum \int (\wtilde{\d}_i \o )^{p/2}  | v_\tau^+ |  \Big) . \label{ff1} \end{align}
Concerning the remainder terms, it follows that
\begin{align}
&  \int  (\wtilde{\d}_i \wtilde{\d}_j)^{p/2} | v_\tau ^+| \leq \| v_\tau^+\| \Big( \int (\wtilde{\d}_i \wtilde{\d}_j )^{n/(n-2)} \Big)^{(n+2)/(2n)} \leq c  \| v_\tau^+\|  \e_{ij}^{\frac{n+2}{2(n-2)}} \ln ^{\frac{n+2}{2n}} (\e_{ij}^{-1}) ,  \label{ka1} \\
& \int_{\mathbb{S}^n} \wtilde{\d}_i ^{p/2}  | v_\tau^+ |  \leq \| v_\tau^+\| \Big( \int_{\mathbb{S}n} \wtilde{\d}_i ^{n/(n-2)} \Big)^{(n+2)/(2n)} \leq c  \| v_\tau^+\|  \frac{\ln ^{(n+2)/(2n)} (\l_i)}{\l_i ^{ (n+2)/4}}  , \label{ka2}
\end{align}
where we have used the following estimate (see $(E2)$  of \cite{B1}) 
\be\label{za1}\int_{\mathbb{S}^n  } (\wtilde{\d}_k \wtilde{\d}_j)^{n/(n-2)} = \int_{\mathbb{R}^n  } (\d_k \d_j)^{n/(n-2)} \leq c \,  \e_{kj} ^{n/(n-2)}\ln(\e_{kj}^{-1}) \quad \forall \, \, k \neq j.\ee
For the other integrals in \eqref{ff1}, we recall that  $ v_\tau ^+ \in N_+(\o)$ and therefore it follows that $ v_\tau ^+ \in H_0(\o)^\perp$. Hence using \eqref{ualphbeta2}, we derive that
\be\label{emna7} -  \langle u_{\a,\b} , v_\tau^+ \rangle +  \int_{\mathbb{S}^n} K u_{\a,\b}^{p-\tau} v_\tau ^+ = -  \langle u_{\a,\b} , v_\tau^+ \rangle +  \int_{\mathbb{S}^n} K u_{\a,\b}^{p} v_\tau ^+  + O(\tau \, \| v_\tau^+\| ) = O(\tau \, \| v_\tau^+\| ) .\ee
Furthermore, using Lemma \ref{lowerL2}, it holds that
\begin{align} \int_{\mathbb{S}^n} K ( \a_i \wtilde{\d}_i)^{p-\tau} v_\tau^+ = & \frac{\a_i^{p-\tau }c_0^{-\tau} K(a_i)}{\l_i ^{\tau(n-2)/2}}  \int_{\mathbb{S}^n} \wtilde{\d}_i ^{p} v_\tau^+ + O\Big(  \int_{\mathbb{S}^n} | K(x) - K(a_i) |   \wtilde{\d}_i ^{p} | v_\tau^+ |  \nonumber  \\
& + \tau \int \wtilde{\d}_i ^p \ln [2+(\l_i ^2 -1)(1-\cos d(a_i,x))] | v_\tau^+ |     \Big) \nonumber\\
& = \frac{\a_i^{p-\tau} c_0^{-\tau} K(a_i)}{\l_i ^{\tau (n-2)/2}}  \langle \wtilde{\d}_i , v_\tau^+ \rangle  + O\Big(   \| v_\tau^+ \| \Big( \tau + \frac{| \n K(a_i) | }{\l_i} + \frac{1}{\l_i ^2}\Big) \Big).\label{emna12} \end{align}
Therefore we get
\be\label{emna8}  \a_i \langle \wtilde{\d}_i , v_\tau^+ \rangle - \int_{\mathbb{S}^n} K ( \a_i \wtilde{\d}_i)^{p-\tau} v_\tau^+ = \a_i \Big( 1-   \frac{\a_i^{p-1-\tau} K(a_i)}{ c_0^{\tau} \l_i ^{\tau (n-2)/2}}  \Big) \langle \wtilde{\d}_i , v_\tau^+ \rangle  + O\Big(   \| v_\tau^+ \| \Big( \tau + \frac{| \n K(a_i) | }{\l_i} + \frac{1}{\l_i ^2}\Big) \Big) . \ee
Finally, using the fact that $v_\tau \in E_{w,a,\l}^\perp$, as in the computations of \eqref{emna11}, we derive that
$$ \langle   \wtilde{\d}_i , v_\tau^+ \rangle = \langle \wtilde{\d}_i , v_\tau \rangle - \langle \wtilde{\d}_i , v_\tau^0 \rangle - \langle \wtilde{\d}_i , v_\tau^- \rangle = O\Big( \| v_\tau \| / \l_i ^{(n-2)/2} \Big) \quad \forall \, \, 1 \leq i \leq N .$$
Thus, using \eqref{qw+-},  \eqref{emna9} becomes
\be\label{emna133} c \| v_\tau^+\| ^2 \leq Q_{\o,a,\l}(v_\tau^+) \leq o( \| v_\tau \| ^2) + O\Big( \| v_\tau \| \Big( R(\tau,a,\l)  \Big) \Big).\ee
Thus using \eqref{emna13},  \eqref{emna5} and   \eqref{emna133}, we derive that
$$ \| v_\tau \| ^2 = \| v_\tau ^+\|^2 + \| v_\tau ^-\|^2 + o( \| v_\tau \|^2) \leq c  \| v_\tau \| \Big(R(\tau ,a,\l)  \Big)$$
which completes the proof.
\end{pf}

Next we are going to give some expansions of the gradient of $I_{K,\tau}$ in $V(\o, N , \mu)$.

\section{Expansion of the gradient in the neighborhood at infinity}

In this section we perform refined asymptotic expansions of the gradient of the Euler-Lagrange functional $I_{K,\tau}$ in the neighborhood at infinity  $V(\o, N , \mu)$. These asymptotic expansions give rise to \emph{balancing conditions} for the parameters of the concentration. We first provide the expansions with respect to the  \emph{gluing parameter } $\a_i$'s. Namely we prove:

\begin{pro}\label{p:35w} Let $u:= u_{\a,\b} + \sum_{ i = 1}^N \a_i  \wtilde{\d}_i + v \in V(\o, N ,\mu )$ with $v$ satisfies \eqref{eq:V0}  and $\tau \ln\l_j$ is small for each $j$.
For  $1\leq i \leq N$, it holds:
$$  \langle \n I_{K,\tau} (u) , \tilde\d_i \rangle  = \a_i S_n \left(1- \l_i^{-\tau \frac{n-2}{2}} \a _i^{\frac{4}{n-2}}K(a_i) \right) + O\Big( R_{\a_i} + \| v \|^2 \Big) $$
where $S_n$ is defined in Theorem \ref{th:blowupConv} and  $R_{\a_i}:= \tau +  \frac{ 1 }{\l_i^2 }+ \sum_{j\neq i} \e_{ij} $ .
\end{pro}

\begin{pf}
We first observe that:
\be \label{nablaI}\langle \n I_{K,\tau} (u) , h \rangle  = \langle u , h \rangle  - \int_{\mathbb{S}^n } K | u | ^{ p-1-\tau} u h \qquad \mbox{ for each } h \in H^1(\mathbb{S}^n ).\ee
 We take $h= \wtilde{\d}_i$ in \eqref{nablaI}. Observe that, since ${v} \in E_{w,a,\l}^\perp$, we get $ \langle {v} , \wtilde{ \d}_i \rangle = 0$. Furthermore, using  the fact that $ u_{\a,\b}$ is $L^\infty$ bounded, easy computations imply that
\begin{align}
& \langle  \wtilde{\d}_j ,  \wtilde{\d}_i \rangle = \int_{\R^n} \d_j ^p \d_i = O( \e_{ij} ) \quad \forall \, \,  j\neq i,  \quad  \langle u_{\a,\b} ,  \wtilde{\d}_i \rangle = \int_{\mathbb{S}^n}  \wtilde{\d}_i ^p ( u_{\a,\b} ) = O\Big( \frac{1}{\l_i^{(n-2)/2}}\Big), \label{az1}\\
& \|  \wtilde{\d}_i \| ^2 = \int_{ \mathbb{S}^n }  \wtilde{ \d }_i ^{p+1} = \int_{\R^n} \d_i ^{p+1} = c_0^{p+1} \int_{\R^n} \frac{ dx}{ (1+| x |^2)^n} := S_n. \label{az2}
\end{align}
For the second part, let $\ov{u}:= u - v$ then it holds
\be\label{az370}  \int_{\mathbb{S}^n} K | u | ^{ p - 1 - \tau} u  \wtilde{\d}_i = \int_{\mathbb{S}^n} K  \ov{u}  ^{ p  - \tau}   \wtilde{\d}_i + ( p - \tau ) \int_{\mathbb{S}^n} K \ov{ u } ^{ p - 1 - \tau} v  \wtilde{\d}_i + O( \| v \|^2 ). \ee
Concerning the first integral, it holds
\be\label{az3}  \int_{\mathbb{S}^n} K \ov{ u } ^{ p  - \tau}   \wtilde{\d}_i =  \a_i^{p-\tau}  \int_{\mathbb{S}^n } K  \wtilde{\d}_i ^{p+1-\tau}  + O\Big(\sum_{j \neq i}  \int_{\mathbb{S}^n } (( \wtilde{\d}_j^p+ u_{\a,\b}^p )  \wtilde{\d}_i +  \wtilde{\d}_i^p ( \wtilde{\d}_j + u_{\a,\b} ) ) \Big). \ee
Note that
\be\label{ww11}
 \int_{\mathbb{S}^n}  \wtilde{\d}_j ^p  \wtilde{\d}_i +  \wtilde{\d}_j  \wtilde{\d}_i ^p = O(\e_{ij}) \quad ;   \quad
 \int_{\mathbb{S}^n} u_{\a,\b}  \wtilde{\d}_i ^p \leq \frac{c} {\l_i^{(n-2)/2}} \quad  ; \quad  \int_{\mathbb{S}^n} u_{\a,\b} ^p  \wtilde{\d}_i \leq c \int_{ \mathbb{S}^n  }  \wtilde{\d}_i  \leq \frac{ c}{\l_i^{(n-2)/2}}.
\ee
Now, since we assumed that $\tau \ln \l_i $ is small then expanding ${K}$ around $a_i$, we get
\begin{align} \int_{\mathbb{S}^n } K  \wtilde{\d}_i ^{p+1-\tau} & =   c_0^{-\tau}  \l_i^{-\tau \frac{n-2}{2}}\int_{ \mathbb{S}^n } K  \wtilde{\d}_i ^{p+1}  + O \Big(  \int _{\mathbb{S}^n } |  \wtilde{\d}_i^{-\tau } - c_0^{-\tau}  \l_i^{-\tau \frac{n-2}{2}} |  \wtilde{ \d}_i ^{p+1} \Big)\label{zza1}\\
& =  c_0^{-\tau}  \l_i^{-\tau \frac{n-2}{2}}\int_{ \mathbb{R}^n } K  {\d}_i ^{p+1} + O \big(  \tau \big) \nonumber \\
& = K(a_i) \l_i^{-\tau (n-2)/2} S_n + O\Big( \tau + \frac{1}{\l_i ^2}  \Big)\nonumber\end{align}
by using \eqref{az2}, Lemma \ref{lowerL2} and the fact that $c_0^{-\tau }  = 1 + O(\tau )$. Hence the estimate of the first integral is complete.\\
For the second integral in the right hand side of \eqref{az370}, using \eqref{est112}, it follows that
\begin{align*}
 \int_{\mathbb{S}^n} K \ov{ u } ^{ p - 1 - \tau} v  \wtilde{\d}_i & = \a_i ^{p-1-\tau} \int_{\mathbb{S}^n} K \wtilde{\d}_i ^{ p  - \tau} v  + O\Big( \sum_{j \neq i} \int_{\mathbb{S}^n} | v | ( \wtilde{\d}_i \wtilde{\d}_j )^{p/2} +  \int_{\mathbb{S}^n} | v | ( \wtilde{\d}_i u_{\a, \b} )^{p/2} \Big)  \\
 & = O\Big(   \| v \| \Big( \tau + \frac{| \n K(a_i) | }{\l_i} + \frac{1}{\l_i ^2} + \sum_{j\neq i} \e_{ij}^{\frac{n+2}{2(n-2)}}\ln^{\frac{n+2}{2n}} (\e_{ij}^{-1})\Big) \Big)
 \end{align*}
by using \eqref{ka1}, \eqref{ka2} and \eqref{emna12} with $v$ instead of $ v_\tau ^+$. Thus the proof follows.
\end{pf}

Next we provide a balancing condition involving the rate of the concentration $\l_i$ and the  mutual interaction of bubbles $\e_{ij}$.
\begin{pro}\label{p:33w}  Let $u:= u_{\a,\b} + \sum_{ i = 1}^N \a_i  \wtilde{\d}_i + v \in V(\o, N ,\mu )$ with $v$ satisfies \eqref{eq:V0} and $\tau \ln\l_j$ is small for each $j$.
For  $1\leq i \leq N$, it holds:
\begin{align*}
& \langle \n I_{K,\tau} (u) , \l_i \frac{\partial\tilde\d_i}{\partial\l_i} \rangle \\
& =  {c_2} \sum _{j\neq i} \a _j \l _i\frac{\partial \e _{ij}}{\partial \l _i} \Big( 1 - \frac{ \a_j ^{p-1} K(a_j) }{\l_j ^{\tau (n-2)/2}} - \frac{ \a_i ^{p-1} K(a_i) }{\l_i ^{\tau (n-2)/2}} \Big)
  + \frac{\a_i ^p}{ \l_i ^{\tau(n-2)/2}}  \Big(c_4 \frac{\Delta K(a_i)}{\l_i^2} + 2c_5 K(a_i) \, \tau \Big) \\
& \quad +O \Big(  \sum_{j\neq k}  \e _{kj }^{\frac{n}{n-2}}\ln (\e _{kj }^{-1}) + \tau ^2 + \sum \frac{\ln\l_k}{\l_k ^{n/2}} +  \frac { | \n K (a_i) |^2 }{\l_i^2} + \frac{1}{\l_i ^3} + \frac{1}{\l_i ^{(n-2)/2}} + \| v \|^2
 \Big)  \end{align*}
 where
 \begin{align*}
  &  c_2:= c_0^{p+1}\int_{\R^n}\frac{1}{(1+|x|^2)^{(n+2)/2}}dx ; \, \, \,   , \\
&  c_4:= \frac{n-2}{2n} c_0^{p+1} \int_{\R^n} \frac{ | x |^2( |x|^2 - 1 )}{(1+| x |^2)^{n+1}} ; \, \, \, \, \,   c_5 := \frac{(n-2)^2}{4} c_0^{p+1}\int_{\R^n}\frac{(|x|^2-1) \ln(1+|x|^2)}{(1+|x|^2)^{n+1}}dx .\end{align*}
\end{pro}

\begin{pf}
We take $h= \l_i \partial \wtilde{\d}_i / \partial \l_i$ in \eqref{nablaI} and we will estimate each term. Observe that, since ${v} \in E_{w,a,\l}^\perp$, it holds: $ \langle \ov{v} , \l_i \partial \wtilde{\d}_i / \partial \l_i \rangle = 0$. Furthermore,
\begin{align}&  \langle \wtilde{\d}_i , \l_i \frac{\partial \wtilde{\d}_i }{ \partial \l_i}  \rangle = \int_{ \R^n} \d_i ^p \l_i \frac{\partial \d_i }{ \partial \l_i}  = 0 , \label{qas11}  \\
& \langle u_{\a,\b} , \l_i \frac{\partial \wtilde{\d}_i }{ \partial \l_i}  \rangle = p \int_{\mathbb{S}^n} \wtilde{\d}_i^{p-1} \l_i \frac{\partial \wtilde{\d}_i }{ \partial \l_i} u_{\a,\b}
= O\Big( \int_{\R^n} \d_i^{p} \Big)  = O\Big(  \frac{1}{\l_i^{(n-2)/2}} \Big). \nonumber
\end{align}
In addition, we have (see Estimate $(E1)$  of \cite{B1}) 
\be \label{v543} \langle \wtilde{\d}_j , \l_i \frac{\partial \wtilde{\d}_i }{ \partial \l_i}  \rangle = c_2 \l_i \frac{\partial \e_{ij}}{\partial \l_i} + O\Big( \e_{ij}^{n/(n-2)} \ln \e_{ij}^{-1} \Big) \quad \forall \, \, j\neq i . \ee
Concerning the second part, let $\ov{u}:= u- {v}$, we have
\be\label{az5}  \int_{\mathbb{S}^n } K | u | ^{ p-1-\tau} u  \l_i \frac{\partial \wtilde{\d}_i }{ \partial \l_i} =  \int_{\mathbb{S}^n } K \ov{u}  ^{ p-\tau}   \l_i \frac{\partial \wtilde{\d}_i }{ \partial \l_i} + (p-\tau) \int_{\mathbb{S}^n } K \ov{u}  ^{ p-1-\tau} {v}  \l_i \frac{\partial \wtilde{\d}_i }{ \partial \l_i} + O(\| {v} \|^2) .\ee
For the second integral, let $\O_i:= \{ x: \sum_{j\neq i} \a_j \wtilde{\d}_j + u_{\a,\b} \leq \a_i \wtilde{\d}_i /2\}$,  it holds
\begin{align}
 \int_{\mathbb{S}^n } K \ov{u}  ^{ p-1-\tau} {v}  \l_i \frac{\partial \wtilde{\d}_i }{ \partial \l_i}  = &  \int_{\mathbb{S}^n } K (\a_i\wtilde{\d}_i)  ^{ p-1-\tau} {v}  \l_i \frac{ \partial \wtilde{\d}_i }{ \partial \l_i} + O\Big(\sum_{j\neq i} \int_{\O_i }\wtilde{\d}_i^{p-1} \wtilde{\d}_j | {v} | \Big)  \label{az6}\\
 & + O\Big(  \int_{\mathbb{S}^n \setminus \O_i }\wtilde{\d}_j^{p-1} \wtilde{\d}_i | {v} | + \int \wtilde{\d}_i^{p-1} u_{\a,\b} | {v} | + \int \wtilde{\d}_i u_{\a,\b} ^{p-1} | {v} | \Big).\nonumber
   \end{align}
 For the last remainder term, using the Holder's inequalities, we obtain
\be\label{vwdelta}  \int \wtilde{\d}_i^{p-1} u_{\a,\b} | v | + \int \wtilde{\d}_i  u_{\a,\b} ^{p-1} | v | \leq c \| v \|  \frac{1 }{ \l _i^2 } . \ee
Observe that, since $n\geq 7$, it follows that $p-1 < 1$. Thus, using the Holder's inequalities and \eqref{za1}, we obtain 
$$\sum_{j\neq i} \int_{\O_i }\wtilde{\d}_i^{p-1} \wtilde{\d}_j | {v} | + \int_{\mathbb{S}^n \setminus \O_i }\wtilde{\d}_j^{p-1} \wtilde{\d}_i | {v} |\leq c \sum_{j\neq i}
\int (\wtilde{\d}_i \wtilde{\d}_j)^{p/2} | {v} | \leq c  \| {v} \| \e_{ij}^{p/2} \ln(\e_{ij}^{-1} )^{\frac{n+2}{2n}} .$$
 Concerning the first integral of the right hand side of \eqref{az6}, following the same computations in \eqref{emna12} and using the fact that $ \langle \partial \wtilde{\d}_i/\partial \l_i , v \rangle = 0$, we obtain
\be  \label{az7} \int_{\mathbb{S}^n } K \d_i ^{p-1-\tau } {v} \l_i \frac{\partial \d_i }{ \partial \l_i}  = O\Big(  \| {v} \| ( \tau  + \frac{|\n K(a_i) |  }{\l_i}  + \frac{1}{\l_i ^2} ) \Big) \ee
which completes the estimate of \eqref{az6}. \\
It remains to estimate the first integral of the right hand side of \eqref{az5}.   Using \eqref{est11}, it holds
\begin{align}
  \int_{\mathbb{S}^n }  K (\ov{u}) ^{p-\tau} \l_i \frac{\partial \wtilde{\d}_i }{ \partial \l_i}  = & \sum_{k=1}^N  \int_{\mathbb{S}^n } K (\a_k \wtilde{\d}_k) ^{p-\tau} \l_i \frac{\partial \wtilde{\d}_i }{ \partial \l_i} +  \int_{\mathbb{S}^n } K (u_{\a,\b}) ^{p-\tau} \l_i \frac{\partial \wtilde{\d}_i }{ \partial \l_i}\label{wdw2} \\
 & + (p-\tau)  \int_{\mathbb{S}^n } K (\a_i \wtilde{\d}_i) ^{p-\tau -1}(\sum_{j \neq i} \a_j \wtilde{\d}_j + u_{\a,\b}) \l_i \frac{\partial \wtilde{\d}_i }{ \partial \l_i} \nonumber \\
 & + \sum_{k\neq j}O\Big(\int (\wtilde{\d}_k \wtilde{\d}_j)^\frac{n}{n-2} + \int (u_{\a,\b} \wtilde{\d}_k)^\frac{n}{n-2} \Big).\nonumber
 \end{align}
The estimate of the first one of the remainder terms is given in \eqref{za1}.
Furthermore, it holds
\be\label{deltaw} \int_{\mathbb{S}^n  }  \wtilde{\d}_{a,\l}^{\frac{n}{n-2}} u_{\a,\b}^{\frac{n}{n-2}}  \leq c  \int_{\mathbb{S}^n }  \wtilde{\d}_{a,\l}^{\frac{n}{n-2}}    \leq c \, \frac{\ln\l }{\l^{n/2}} . \ee
In addition, using \eqref{za1}, \eqref{v543} and Lemma \ref{lowerL2},  we get, for $k \neq i$,
\begin{align*}  \int_{\mathbb{S}^n } K  &  \wtilde{\d}_k ^{p-\tau} \l_i \frac{\partial \wtilde{\d}_i }{ \partial \l_i} \\
&  = \frac{c_0^{-\tau}}{\l_k^{\tau(n-2)/2}} K (a_k)  \int_{\mathbb{S}^n }  \wtilde{\d} _k ^{p} \l_i \frac{\partial \wtilde{\d}_i }{ \partial \l_i}  + O\Big(  \int_{\mathbb{S}^n } d(x,a_k)  \wtilde{\d}_k ^{p}   \wtilde{\d}_i  + \int_{\mathbb{S}^n } | \wtilde{\d}_k ^{-\tau } - \frac{c_0^{-\tau}}{\l_k ^{\tau (n-2)/2}} |  \wtilde{\d}_k ^{p}  \wtilde{\d}_i \Big)\\
& = \frac{1}{\l_k^{\tau(n-2)/2}} K (a_k)  c_2 \l_i \frac{\partial \e_{ik}} {\partial \l_i} + O\Big(  \e_{ik}^{n/(n-2)}\ln\e_{ik}^{-1} +( \frac{1}{\l_k} + \tau ) \, \e_{ik}  (\ln\e_{ik}^{-1} )^{(n-2)/n}  \Big)
 .\end{align*}
In the same way, it holds
$$  p \int_{\mathbb{S}^n } K  \wtilde{\d}_i ^{p-\tau-1}  \wtilde{\d}_k  \l_i \frac{\partial \wtilde{\d}_i }{ \partial \l_i} =\frac{1}{\l_i ^ {\tau \frac{n-2}{2}}} K (a_i)  c_2 \l_i \frac{\partial \e_{ik}} {\partial \l_i} + O\Big(  \e_{ik}^\frac{n}{n-2}\ln\e_{ik}^{-1} +( \frac{1}{\l_i} + \tau ) \, \e_{ik} ( \ln\e_{ik}^{-1} )^\frac{n-2}{n}  \Big).$$
 Concerning the first term of \eqref{wdw2} for $k =i$,  using Lemma \ref{lowerL2}, we get
\begin{align}   \int_{\mathbb{S}^n } K  \wtilde{\d}_i ^{ p-\tau }   \l_i \frac{\partial \wtilde{\d}_i }{ \partial \l_i}  = & \frac{ c_0^{-\tau } }{\l_i^{\tau \frac{n-2}{2}}} \Big(  \int_{\mathbb{S}^n } K  \wtilde{\d}_i  ^{ p}   \l_i \frac{\partial \wtilde{\d}_i }{ \partial \l_i}  \nonumber \\
 & + \tau \frac{n-2}{2}  \int_{\mathbb{S}^n } K  \wtilde{\d}_i  ^{ p}   \l_i \frac{\partial \wtilde{\d}_i }{ \partial \l_i}\ln[2+ (\l_i^2-1)(1-\cos d(  x,a_i ))] \Big) + O( \tau ^2 ). \label{dd1}
\end{align}
In addition, using \eqref{qas11}, it holds
\begin{align*}
 \int_{\mathbb{S}^n } K  \wtilde{\d}_i  ^{ p}   \l_i \frac{\partial \wtilde{\d}_i }{ \partial \l_i}  & = \int_{\R^n } \tilde{K}  \d_i  ^{ p}   \l_i \frac{\partial \d_i }{ \partial \l_i} =
 \frac{1}{2} \sum_{1\leq k,j \leq n} \frac{\partial ^2 \tilde{K}}{ \partial x_k \partial x_j }(\tilde {a}_i) \int_{\R^n } (x-\tilde{a}_i)_k   (x-\tilde{a}_i)_j  \d_i  ^{ p}   \l_i \frac{\partial \d_i }{ \partial \l_i} + O(\frac{1}{\l_i^3})\\
&  =   \frac{1}{2}   \sum_{1\leq k \leq n} \frac{\partial ^2 \tilde{K}}{ \partial x_k ^2}(\tilde {a}_i) \int_{\R^n } (x-\tilde{a}_i)_k ^2    \d_i  ^{ p}   \l_i \frac{\partial \d_i }{ \partial \l_i} + O(\frac{1}{\l_i^3})\\
& =  \frac{1}{2}  \D  \tilde{K} (\tilde {a}_i)   \frac{n-2}{2} \frac{c_0^{2n/(n-2)}}{\l_i ^2 }  \frac{1}{n}   \int_{\R^n} \frac{ | y |^2 (1-| y |^2)}{(1+ | y |^2)^{n+1}}dy + O(\frac{1}{\l_i^3}),  \end{align*}

\begin{align}
& \int_{\mathbb{S}^n } K  \wtilde{\d}_i  ^{ p}   \l_i \frac{\partial \wtilde{\d}_i }{ \partial \l_i}\ln [ 2+ (\l_i^2-1)(1-\cos d ( x,a_i )) ]  =   \int_{\mathbb{R}^n } \tilde{K}  \d_{0,\l_i}  ^{ p}   \l_i \frac{\partial \d_{0,\l_i} }{ \partial \l_i}\ln [ \frac{(2+2\l_i^2|y|^2)}{1+|y|^2} ] \label{rem1}  \\
& \qquad \qquad \qquad  = \int_{\mathbb{R}^n } \tilde{K}  \d_{0,\l_i}  ^{ p}   \l_i \frac{\partial \d_{0,\l_i} }{ \partial \l_i}\ln (2+2\l_i^2|y|^2) + O \Big( \int_{\mathbb{R}^n } \d_{0,\l_i}  ^{ p+1}  \ln (1+|y|^3) \Big)  \nonumber \\
& \qquad \qquad \qquad  = K(a_i) c_0^\frac{2n}{n-2} \frac{n-2}{2}  \int_{\R^n} \frac{ (1-| y |^2)\ln(1+ | y |^2)}{(1+ | y |^2)^{n+1}}dy +O(\frac{1}{\l_i ^2}), \nonumber
\end{align}
which achieve  the estimate of \eqref{dd1}. Note that, for  the estimate of the remainder term in \eqref{rem1} we used the fact that $ \ln ( 1 + | t | ^3 ) \leq c | t |^3 $ for each $ t \in \R$.\\
It remains to estimate the integrals which appear in \eqref{wdw2} and  contain $u_{\a,\b}$. Using \eqref{ww11}, there hold
$$ \int_{\mathbb{S}^n } K  ( \a_0 u_{\a,\b} )  ^{ p-\tau }   \l_i \frac{\partial \wtilde{\d}_i }{ \partial \l_i} =  O\Big(  \int _{\mathbb{S}^n }   \wtilde{\d}_i \Big) =  O\Big( \frac{ 1 }{ \l_i ^{(n-2)/2}} \Big), $$
\begin{align*}  \int_{\mathbb{S}^n } K  & \wtilde{\d}_i ^{ p-\tau -1}  u_{\a,\b}    \l_i \frac{\partial \wtilde{\d}_i }{ \partial \l_i}
 = O\Big(  \int_{\mathbb{S}^n }    \wtilde{\d}_i ^{ p}  \Big)   =   O\Big( \frac{ 1 }{ \l_i ^{(n-2)/2}} \Big).
\end{align*}
Combining the previous estimates, the proof follows.
\end{pf}

Finally, we provide the following balancing condition involving  the point of concentration $a_i$.
  \begin{pro}\label{p:38w}   Let $u:= u_{\a,\b} + \sum_{ i = 1}^N \a_i  \wtilde{\d}_i + v \in V(\o, N ,\mu )$ with $v$ satisfies \eqref{eq:V0}  and $\tau \ln\l_j$ is small for each $j$.
For  $1\leq i \leq N$, it holds:
 $$
\langle \n I_{K,\tau} (u) , \frac{1}{\l_i}  \frac{\partial\tilde\d_i}{\partial a_i} \rangle =  \a_i^p ( c + o(1) ) \frac{\nabla K(a_i)}{\l_i} +O \Big(\frac{1}{\l_i^3} + \frac{1 }{\l_i ^{(n-2)/2}} + \sum_{j\neq i}\e_{ij} +  \| v \|^2 \Big) . $$
\end{pro}
 \begin{pf}
The proof follows as the previous one. Hence we will omit it here.
\end{pf}

In the next section,  we are going to prove Theorem \ref{th:blowup}.

 \section{Asymptotic behavior of blowing up solutions}

The focus in this section is a refined blow up analysis of finite energy  solutions of the  subcritical approximated problem $(\mathcal{N}_{K,\tau})$. Namely, our aim is to prove Theorem \ref{th:blowup}.  To do this, we need to fix some important facts. \\
 Let $(u_\tau)$ be a sequence of energy bounded solutions of $(\mathcal{N}_{K,\tau})$ converging weakly but not strongly $u_\tau \rightharpoonup \o \neq 0$. It follows from Propositions \ref{p:35w}, \ref{p:33w} and \ref{p:38w} (by taking $u_\tau$ instead of $u$) that the following system holds:
\begin{align}
&| 1- \l_i^{-\tau \frac{n-2}{2}} \a _i^{\frac{4}{n-2}}K(a_i) | = O\Big( R_{\a_i} + R(\tau,a,\l)^2 \Big)\quad \forall \, \, i \leq N, \nonumber\\
 & (E_i) \qquad - {c_2} \sum _{j\neq i} \a _j \l _i\frac{\partial \e _{ij}}{\partial \l _i}    +  {\a_i} \Big(c_4 \frac{\Delta K(a_i)}{\l_i^2 K(a_i)} + 2c_5 \tau \Big) = O(R_{2,\l}) \quad \forall \, \, i \leq N, \nonumber\\
& (F_i) \qquad \frac{|\nabla K(a_i)|}{\l_i} \leq c \Big(\frac{1}{\l_i^3} + \sum_{j\neq i}\e_{ij} + R^2(\tau,a,\l) + \sum \frac{\ln \l_k }{\l_k ^{n/2}} \Big) \quad \forall \, \, i \leq N , \nonumber
  \end{align} 
  where $R(\tau,a,\l)$ is defined in Proposition \ref{eovvw} and 
  \be \label{R2l}  R_{2,\l} :=   \sum_{j\neq k}  \e _{kj }^{\frac{n}{n-2}}\ln (\e _{kj }^{-1}) + \tau ^2 + \sum \frac{1}{\l_k ^{5/2}} + \sum \frac { | \n K (a_k) |^2 }{\l_k ^2} . \ee

 Now, we order the $\l_i$'s as: $ \l_1 \leq \l_2 \leq \cdots\leq \l_N$ and we start by the following estimate.

 \begin{lem}\label{sum} For each solution $u_\tau$ in $V(\o, N,\tau)$, it holds
 $$\sum_{k\neq i} \e_{ik} + \sum  \frac{ | \n K(a_k) | }{ \l_k}  + \tau \leq c  \frac{1}{\l^2_1}.$$
 Furthermore, $R_{2,\l}$, defined in \eqref{R2l}, satisfies : $R_{2,\l} =  o(1/\l_1^{2(n-1)/(n-2)})$.
 \end{lem}
 \begin{pf}
 We first notice that
 \begin{align}\label{derepsilon}
& - \l_i \frac{\partial \e_{ij}}{\partial \l_i} - \l_j \frac{\partial \e_{ij}}{\partial \l_j
} \geq 0 \, \, \mbox{ for each } i\neq j  \quad \mbox{ and } \quad  - \l_i \frac{\partial \e_{ij}}{\partial \l_i} \geq c \e_{ij} \quad \mbox{ if } \l_i \geq c' \l_j.\end{align}
Hence, summing $2^i \a_i(E_i)$, we obtain
$$  \sum_{ k\neq i} \e_{ik} +   \tau  \leq c   \Big( R_{2,\l} + \frac{1}{\l^2_1}\Big) .$$
The result follows by using $(F_i)$  and the definition of $R_{2,\l}$ given in \eqref{R2l}.
\end{pf}

 \begin{lem}\label{56}
 For $i\in \{1, \cdots, N\}$, we define  \be\label{Gi}\G_i:= \l_i^2 \sum_{k \neq i} \e_{ik}  + \frac{| \n K(a_i) | /\l_i }{\sum \e_{ki} + 1/\l_i^2}.\ee
 Let $D_1':=\{ i : \lim \G_i =\infty\}$ and $D_1: = \{1,\cdots,N\} \setminus D_1'$. \\
$(i)$  For each $i \in D_1$, there exists a critical point  $y_i$ of $K$ such that $\l_i | a_i - y_i | \leq C$.\\
$(ii)$ If there exists an index $i \in D_1'$, then it holds that
$$ \sum_{ j \geq i} \Big( \frac{ | \n K(a_i) | }{\l_i} + \sum_{k\neq j} \e_{kj} + \tau   \Big) = o\Big(\frac{1}{\l_1^{2(n-1)/(n-2)}}\Big). $$
$(iii)$ There exists a critical point $ y $ of $ K $ such that $ a_1 $ converges to $y$ .
 \end{lem}
 \begin{pf}
 Let $i\in D_1$, it follows that $\G_i$ is bounded and therefore it holds that $| \n K(a_i) | \leq C/\l_i$ which implies the first assertion. Concerning the second one, let $i\in D_1'$, summing $2^j \a_j (E_j)+(F_j)/m$ for $j \geq i$  where $m$ is a small constant, it holds that
 $$ \sum_{ j \geq i} \Big( \frac{ | \n K(a_j) | }{\l_j} + \sum_{k\neq j} \e_{kj} + \tau   \Big) = O\Big( \frac{1}{\l_i^2} + R_{2,\l} \Big).$$
 Since $i\in D_1'$, it follows that $1/\l_i ^2$ is small with respect to the left hand side. Thus the proof follows from the estimate of $R_{2,\l}$ (see Lemma \ref{sum}). \\
 Concerning the proof of Claim $(iii)$, observe that, if $ 1 \in D_1 $, the result follows from Claim $(i)$, however, if $ 1 \in D_1' $, from Claim $(ii)$, we deduce that $ | \n K(a_1) |  = o( \l_1^{-n/(n-2) } ) $ which implies the result in this case. Thus the proof of the lemma is complete. 
 \end{pf}

 \begin{lem}\label{55}
  For each $i \neq  j $, it holds that : $\e_{ij}= o(1/\l^2 _1)$.
 \end{lem}
 \begin{pf}
  Observe that, if $i$ or $j$ belongs to $D_1'$, then the result follows from Lemma \ref{56}. In the other case, that is $i,j \in D_1$, using again Lemma \ref{56}, there exist critical points $y_i$ and $y_j$ such that $\l_k d( a_k, y_k ) \leq c $ for $k=i,j$. Two cases may occur:
\begin{enumerate}
\item  [(a)] Either $y_i \neq y_j$, and in this case we get $d( a_i, a_j ) \geq c $ and therefore the result follows easily,
\item [(b)] or $y_i = y_j$. Since we have $\l_k  d( a_k, y_k ) \leq c $ and $\e_{ij}$ is small, it follows that $\l_i / \l_j \to 0$ or $\infty$. Taking $\l_i \leq \l_j$ and using the fact that $\G_j$ is bounded, it holds
 $$ \e_{ij} \leq \frac{c}{  \l_j ^2}  = c \frac{ \l_i ^2}{  \l_j ^2} \frac{1}{  \l_i ^2} = o \big( \frac{1}{  \l_i ^2} \big) = o \big( \frac{1}{  \l^2_1} \big).$$
 \end{enumerate}
 Hence  the proof is completed.
 \end{pf}

 Next we are going to complete the proof of Theorem \ref{th:blowup}.  Using the above lemmas, Equations $(E_i)$  can be improved and we get
 \be\label{Ei'} (E_i') :  \qquad   2c_5 \tau + c_4 \frac{ \D K(a_i) }{\l_i ^2 K(a_i)} = o\big( \frac{1}{\l^2 _1} \big)    . \ee

Now, we claim that
\be\label{claim} (a) \quad \frac{\l_N}{\l_1} \leq C \qquad \quad \mbox{ and } \qquad  \quad (b) \quad  D_1' = \emptyset .\ee
Arguing by contradiction, assume that $ \l_N / \l_1 \to \infty$. It follows from $(E_{N}')$ that $ \tau = o( 1/\l_1^2)$.   Putting this information in $(E_1')$ and using Claim $(iii)$ of Lemma \ref{56}, we derive a contradiction. Hence all the $ \l_i$'s are of the same order.
Concerning assertion $(b)$, it follows from assertion $(a)$ and Lemma \ref{sum}.  This completes the proof of the claim. 

Combining \eqref{claim} and assertion $(i)$ of Lemma \ref{56}, we derive that, for each $i$, there exists $y_{j_i}$ such that $\l_i d(a_i , y_{j_i}) \leq c$. Assume that, for $i \neq k$, we have $y_{j_i} = y_{j_k}$. Thus we obtain that $\l_i \l_k d(a_i,a_k)^2 \leq c$. This information and the fact that the $ \l_i$'s are of the same order contradict the smallness of $\e_{ik}$. Therefore we derive that $ d(a_i,a_k) \geq c $ for each $i \neq k$ which implies that $ \e_{ik} = O (1/\l_i ^{n-2})$ for each $k\neq i$. Now, putting the previous informations in the right hand side of $(F_i)$ we derive that $ \l_i d( a_i , y_{j_i} ) \to 0$ as $\tau \to 0$.

Regarding \eqref{lll}, it follows from $(E_i')$ and the fact that the $ \l_i$'s are of the same order.

Finally, concerning the morse index of the solution $u_\tau$, the result is done for $\o = 0$ in \cite{MM}. However, for $\o \neq 0$ and non-degenerate, it follows that $u_{\a,\b} = \a \o$.  In this case, we need to add the contribution of $\a$ and $\o$ in the morse index formula. Note that, by easy computations, we see that
$$ 
D^2 I_{K,\tau} (u_\tau) (\frac{\partial u_\tau}{\partial \a}, \frac{\partial u_\tau}{\partial \a}) = D^2 I_{K,\tau} (u_\tau) ( \o , \o )+ O\left( \Big| \Big|\frac{\partial v_\tau}{\partial \alpha}\Big|\Big|\right).
$$
Note that $D^2 I_{K,\tau} (u_\tau) ( \o , \o ) \leq -c<0$ and $\Big| \Big|\frac{\partial v_\tau}{\partial \alpha}\Big|\Big| = o(1)$ (see Lemma \ref{Lemmalast}).
This implies that the variable $\a$ belongs to the negative space. For the contribution of $\o$, using Lemma \ref{Qwalpositive}, we see that it is equal to
$$ \mbox{ dim } (N_-(\o)) = \mbox{ morse index } (I_{K} , \o) .$$
This completes the proof of the theorem.

\section{Construction of blowing up solutions for the approximate problem}

Let $n \geq 7$, $0 < K \in C^3({\mathbb{S}^n})$ and let $y_1, \cdots, y_N$ be distinct non-degenerate critical points of $K$ with $\D K(y_i) <0$ and $\o$ be a solution of $(\mathcal{N}_K)$. In Theorem \ref{th:blowupConv} we assume that $\o$ is non-degenerate and therefore the variable $\beta$ introduced in the function $u_{\a,\b}$ has to be equal to zero, that is, $u_{\a,\b}=\a\o$.

As in \cite{BLR}, the strategy of the proof is the following: we will define a set $M_\tau$ whose elements  are some points $\mathcal{M} :=(\a,\l,x,v) \in (\R_+)^{N+1} \times  (\R_+)^N\times ({ \mathbb{S}}^n)^N \times H^1(\mathbb{S}^n)$  where $v$ satisfies \eqref{eq:V0}, the other variables satisfy some conditions and  $N$ is the number of the bubbles in the desired constructed solution. Precisely, let

\begin{align}  M_\tau :=  \{  \mathcal{M} := &  (\a_0, \a,\l,x,v)\in \R_+ \times (\R_+)^N \times(\R_+)^N\times ( \mathbb{S}^n)^N \times H^1(\mathbb{S}^n):  \label{Meps}\\
& v\mbox{ satisfies \eqref{eq:V0}} ; \, \,   |{\a_i^{4/(n-2)}K(x_i)} - 1 | < \tau \ln^2\tau\,  \, \forall \, \, i \geq 1;  \nonumber\\
&  |\a_0^{4/(n-2)} - 1 | < \tau \ln^2\tau\,   \, ; \, \,   C^{-1} \tau \leq \l_i^{-2} \leq  C \tau \, \,  ; \, \, \,   d(x_i , y_i )  \leq C \tau \,\, \forall \, i\geq 1   \}, \nonumber
\end{align}
where $C$ is a large positive constant.

 We remark that, for  $\mathcal{M} :=   (\a_0, \a,\l,x,v) \in M_\tau$, it holds that $  u:=  \a_0 \o + \sum_{i=1}^{N} \a_i \wtilde{\d}_{ x_i,\l_i } + v \in V(\o, N , \mu)$ for some $\mu $ small.

In addition, we define a function
\begin{eqnarray}\label{F10}  \Psi_{\tau}:M_{\tau}\,\rightarrow \R ;\quad \mathcal{M}=(\a_0, \a,\l,x,v)\mapsto I_{K,\tau} \big( \a_0 \o + \sum_{i=1}^{N} \a_i \wtilde{\d}_{x_i,\l_i} \, + \, v\big) \end{eqnarray}
and we need to find a critical point of $ \Psi_{\tau}$.

Recall that the variable $v$ satisfies some orthogonality conditions. Thus using the Lagrange multiplier theorem, it is easy to get the following proposition.

\begin{pro}\label{G1} Let $ \mathcal{M} =(\a_0, \a,\l,x,v)\in M_{\tau}$. $\mathcal{M}$ is a critical point of $\Psi_{\tau}$ if and only if $ \a_0 \o + \sum_{i=1}^{N}\a_i\wtilde{\d}_{x_i,\l_i}   +  v$ is a critical point of $I_{K,\tau}$, i.e. if and only if there exists $\big(A_0,A,B,C\big)\in \R\times \R^N\times\R^N\times\big(\R^{n}\big)^N$ such that the following holds :
\begin{align}
 & (E_{\a_i}) \qquad \frac{\partial \Psi_{\tau}}{\partial \a_i}(\a_0, \a,\l,x,v)=0,\,\, \forall \,i=0, 1,\cdots,N ,  \notag\\
 & (E_{\l_i}) \qquad \frac{\partial \Psi_{\tau}}{\partial \l_i}(\a_0, \a,\l,x,v)=B_i\langle \l_i \frac{\partial^2 \wtilde{\d}_i}{\partial\l_i^2},v\rangle+\sum_{j=1}^{n} C_{ij}\langle \frac{1}{\l_i}\frac{\partial^2 \wtilde{\d}_i}{\partial x_i^j\partial\l_i},v\rangle,\,\ \forall \,i=1,\cdots, N , \notag\\
 & (E_{x_i}) \qquad \frac{\partial \Psi_{\tau}}{\partial x_i}(\a_0, \a,\l,x,v)\lfloor_{T_{x_i}( \mathbb{S}^n)} \, = \, B_i\langle \l_i \frac{\partial^2 \wtilde{\d}_i}{\partial\l_i\partial x_i},v\rangle+ \sum_{j=1}^{n} C_{ij} \langle \frac{1}{\l_i}\frac{\partial^2 \wtilde{\d}_i}{\partial x_i^j\partial x_i},v\rangle, \,\ \forall \,i=1,\cdots,N ,  \notag\\
 & (E_v)\qquad \frac{\partial \Psi_{\tau}}{\partial v}(\a_0, \a,\l,x,v)=A_0 \o + \sum_{i=1}^N \biggl( A_i \wtilde{\d}_i+B_i\l_i \frac{\partial \wtilde{\d}_i}{\partial \l_i}+\sum_{j=1}^{n} C_{ij}\frac{1}{\l_i}\frac{\partial \wtilde{\d}_i}{\partial x_i^j}\biggr)\notag
\end{align}
where $\wtilde{\d}_i:=\wtilde{\d}_{x_i,\l_i} $.
\end{pro}
The results of Theorem \ref{th:blowupConv} will be obtained through a careful analysis of the previous equations on $M_{\tau}$. Note that
\begin{align*} & \frac{\partial \Psi_{\tau}}{\partial v}(\a_0, \a,\l,x,v) =   \n I_{K,\tau} (u); \\
& \frac{\partial \Psi_{\tau}}{\partial \a_0}(\a_0, \a,\l,x,v) = \langle \n I_{K,\tau}  (u) , \o \rangle \end{align*}
and for each $ i \geq 1$, we have
\begin{align*} & \frac{\partial \Psi_{\tau}}{\partial \a_i}(\a_0, \a,\l,x,v) = \langle \n I_{K,\tau}  (u) , \tilde\d_i \rangle  ; \\
& \frac{\partial \Psi_{\tau}}{\partial x_i}(\a_0, \a,\l,x,v)_{ \lfloor_{T_{x_i}( \mathbb{S}^n)}} \, = \, \langle \n I_{K,\tau}  (u) , \a_i \frac{\partial\tilde\d_i}{\partial x_i} \rangle_{ \lfloor_{T_{x_i}( \mathbb{S}^n)}}  \quad  ; \\
& \frac{\partial \Psi_{\tau}}{\partial \l_i}(\a_0, \a,\l,x,v) \, = \, \langle \n I_{K,\tau}  (u) , \a_i \frac{\partial\tilde\d_i}{\partial \l_i} \rangle \end{align*}
where $ u := \a_0 \o + \sum_{i=1}^{N} \a_i \wtilde{\d}_{x_i,\l_i} \, + \, v.$

We remark that  the last three equations are estimated in Propositions \ref{p:35w}, \ref{p:33w} and \ref{p:38w}. 

\begin{pro}\label{ovvmin}
Let $ u:=  \a_0 \o + \sum_{i=1}^{N} \a_i \wtilde{\d}_{ x_i,\l_i } \in V(\o, N, \mu)$. Then there exists a unique $\ov{v} \in E_{\o, x, \l}^\perp$ satisfying
\be\label{wcrit} \langle \n  I_{K,\tau}  ( u + \ov{v}) , h \rangle = 0 \qquad  \forall \, \, h \in E_{\o, x, \l} ^\perp .\ee
This function $\ov{v}$ satisfies
$$ \| \ov{v} \| \leq c R(\tau ,a,\l) $$ 
where $ R(\tau ,a,\l) $ is defined in Proposition \ref{eovvw}.
\end{pro}

\begin{pf}
To prove \eqref{wcrit}, we need to expand $I_{K,\tau} ( u + {v})$ with respect to $v \in  E_{\o,a,\l}^\perp$. Recall that
$$ I_{K,\tau} (u+v) := \frac{1}{2} \| u+v \|^2 - \frac{1}{p+1 - \tau} \int_{ \mathbb{S}^n} K | u + v | ^{p+1-\tau} .$$
Note that, for each $x,y \in \R$ and $\gamma \in (2,3)$, we have
$$ | x+y|^\g = | x | ^\g + \g | x |^{\g-2} x y + \frac{1}{2} \g ( \g -1 ) | x |^{\g-2} y ^2 + O( | y | ^\g).$$
For $ n \geq 7$ and $\tau$ small, it holds that $ p + 1 - \tau \in (2,3)$ and therefore we get
\begin{align*}  \int_{ \mathbb{S}^n} K | u + v | ^{p+1-\tau}  = &  \int_{ \mathbb{S}^n} K | u  | ^{p+1-\tau} + (p+1-\tau) \int_{ \mathbb{S}^n} K | u  | ^{p-1-\tau} u v \\
& +  \frac{1}{2} (p+1-\tau) (p-\tau)  \int_{ \mathbb{S}^n} K | u  | ^{p-1-\tau}  v^2 + O \Big( \int_{ \mathbb{S}^n} K | v  | ^{p+1-\tau} \Big). \end{align*}
Since $ u \in V(\o, N , \mu )$ with $\tau \ln \l_i $ is small for each $i$, it follows that
$$ \int_{ \mathbb{S}^n} K | u  | ^{p-1-\tau}  v^2 = \sum_{i=1}^N \int_{ \mathbb{S}^n} {\wtilde{\d}_i} ^{p-1}  v^2  + \int_{ \mathbb{S}^n} K  \o^{p-1}  v^2 + o( \| v \| ^2).$$
Hence we get
\be\label{ikt}  I_{K,\tau} (u+v) = I_{K,\tau} (u) - f(v) + \frac{1}{2} Q_{\o, a , \l } (v) + o( \| v \| ^2) \ee
where
$$ f(v) := \int_{ \mathbb{S}^n} K | u  | ^{p-1-\tau} u v \quad \mbox{ and }  \quad
  Q_{\o, a , \l } (v)  := \| v \|^2 - p \sum_{i=1}^N \int_{ \mathbb{S}^n} {\wtilde{\d}_i} ^{p-1}  v^2  - p \int_{ \mathbb{S}^n} K  \o^{p-1}  v^2 .$$
Thus \eqref{wcrit} follows from Proposition \ref{Qwalpositive} (which tells us that $Q_{\o, a , \l }$ is non-degenerate on the space $E_{\o, a , \l }^\perp$). \\
Now it remains to estimate $ \| \ov{v } \|$. Note that, using the Lagrange multiplier theorem, \eqref{wcrit} implies the existence of some constants $(A_0,A,B,C) \in \R \times \R^N \times \R^N \times (\R^n)^N $ such that
\be\label{qas1} \n I_{K,\tau} (u+ \ov{v} ) = A_0 \o + \sum_{i=1}^N \Big( A_i \wtilde{\d}_i + B_i \l_i \frac{\partial  \wtilde{\d}_i}{\partial \l_i} + C_i \cdot  \frac{1}{\l_i} \frac{\partial  \wtilde{\d}_i}{\partial x_i} \Big) \ee
where $ \wtilde{\d}_i := \wtilde{\d}_ {x_i, \l_i} $.
Since $\o$ is assumed to be non-degenerate, we derive that
\be \label{val} H^1(\mathbb{S}^n) = N_-(\o) \oplus \mbox{ span}\{\o \} \oplus N_+( \o ), \ee
that is $H_0(\o)$ (defined in \eqref{vvv*1}) becomes $H_0(\o) = \mbox{ span}\{\o \}.$\\
As in the proof of Proposition \ref{eovvw}, we divide $\ov{v}$ into
$$ \ov{v} = \ov{v}_- + \ov{v}_0 + \ov{v}_+$$
and we have $ \| \ov{v}_0 \| = o( \| \ov{v} \| )$ (by using Proposition \ref{Qwalpositive}).\\
Multiplying \eqref{qas1} by $\ov{v}_-$ and using the fact that $ \langle \ov{v}_- , \o \rangle = 0$, we get
$$ \langle  \n I_{K,\tau} (u+ \ov{v} ) , \ov{v}_- \rangle = \sum O \Big( \int_{\mathbb{S}^n} \wtilde{\d}_i ^p | \ov{v}_- | \Big) =  \sum O \Big( \frac{ \| \ov{v}_- \| } { \l_i ^{(n-2)/2}} \Big)  $$
(since $ \ov{v}_- $ belongs to a fixed finite dimensional space which implies that $ \| \ov{v}_- \|_\infty \leq c \| \ov{v}_ -\| $). \\
The left hand side is computed in \eqref{emna1} with $ v_\tau $ instead of $ \ov{v }$. Following the computations done in \eqref{emna1}-\eqref{emna5}, we derive that
$$ \| \ov{v}_- \| ^2 \leq  c \, \| \ov{v}_- \| \Big( \tau + \sum \frac{1}{\l_i ^{(n-2)/2}} + o( \| \ov{v} \| ) \Big). $$
Now multiplying  \eqref{qas1} by $\ov{v}_+$ and using the fact that $ \langle \ov{v}_+ , \o \rangle = 0$, we get
$$ \langle  \n I_{K,\tau} (u+ \ov{v} ) , \ov{v}_+ \rangle = \sum \Big( A_i \langle  \wtilde{\d}_i  , \ov{v}_+ \rangle + B_i \langle \l_i \frac{\partial  \wtilde{\d}_i}{\partial \l_i} , \ov{v}_+ \rangle +  \langle C_i \cdot  \frac{1}{\l_i} \frac{\partial  \wtilde{\d}_i}{\partial a_i} , \ov{v}_+ \rangle \Big) . $$
The left hand side is computed in \eqref{emna1} with $ v_\tau $ instead of $ \ov{v }$.
Recall that $\ov{v}_+$ does not belong to $ E_{\o, a , \l}^\perp$ but $ \ov{v} \in E_{\o, a , \l}^\perp$. Thus
$$ \langle  \wtilde{\d}_i  , \ov{v}_+ \rangle = \langle  \wtilde{\d}_i  , \ov{v} \rangle - \langle  \wtilde{\d}_i  , \ov{v}_- \rangle - \langle  \wtilde{\d}_i  , \ov{v}_0 \rangle = O \Big( \int \wtilde{\d}_i ^p ( | \ov{v}_- | +  | \ov{v}_0 | ) \Big) =  O \Big( \frac{ \| \ov{v} \| } { \l_i ^{(n-2)/2}} \Big). $$
The same holds for the other scalar products.\\
Following the computations done in the proof of Proposition \ref{eovvw}, we get
$$ \| \ov{v}_+ \| ^2 \leq o( \| \ov{v} \|^2 ) +  c \, \| \ov{v} \| R(\tau, a , \l), $$
where $ R(\tau, a , \l)$ is introduced in Proposition \ref{eovvw}.\\
Finally, since $ \| \ov{ v } \|^2 =  \| \ov{ v }_+ \|^2 +  \| \ov{ v }_- \|^2  +  \| \ov{ v }_0 \|^2 $ we derive the desired estimate. \\
The proof of Proposition \ref{ovvmin} is thereby completed.
\end{pf} 

\begin{pro}\label{p:1conv}
 Let $ u := \a_0 \o + \sum_{i=1}^{N} \a_i \wtilde{\d}_{ x_i,\l_i } \, + \, \ov{v}$. For each $1\leq i \leq N$, It holds:
$$
 \langle \n  I_{K,\tau} (u),\o \rangle =  \a_0\| \o \|^2 \, (1-\a_0^{p-1})+O\Big( \tau + \| \ov{v} \|^2 + \sum \frac{1}{\l_i ^{(n-2)/2}}+ \frac{1}{\l_i ^4} \Big) .$$
\end{pro}
\begin{pf}
Since $\ov{v}$ satisfies \eqref{eq:V0}, using \eqref{nablaI} and \eqref{az1}, it holds
$$ \langle \n   I_{K,\tau} (u),\o \rangle = \langle u,\o \rangle - \int K | u |^{p-\tau-1} u \o = \a_0 \| \o \| ^2 + \sum O \Big( \frac{1}{\l_i ^{(n-2)/2}} \Big) - \int K | u |^{p-\tau -1} u \o .$$
Furthermore, let $ \ov{u} := \a_0 \o + \sum_{i=1}^{N} \a_i \wtilde{\d}_{ x_i,\l_i } $, we have $$\int_{\mathbb{S}^n } K | u | ^{ p-1-\tau} u  \o =  \int_{\mathbb{S}^n } K \ov{u}  ^{ p-\tau}   \o + (p-\tau) \int_{\mathbb{S}^n } K \ov{u}  ^{ p-1-\tau} \ov{v}  \o  + O(\| \ov{v} \|^2) . $$
Using \eqref{vwdelta} and the fact that $ \langle w ,  \ov{v} \rangle = 0 $, it follows that 
\begin{align*} \int_{\mathbb{S}^n } K \ov{u}  ^{ p-1-\tau} \ov{v}  \o & = \a_0^{p-1-\tau} \int_{\mathbb{S}^n } K \o ^{p-\tau} \ov{v} + \sum_{j\geq 1} O\Big( \int \o  ^{p-1} \wtilde{\d}_j | \ov{v} | + \int \o  \wtilde{ \d }_j ^{p-1} | \ov{v} | \Big) \\
 & = O \Big( \| \ov{ v } \|  ( \tau + \sum \frac{ 1 } { \l_j ^2 } )  \Big) . 
 \end{align*}
In addition, using \eqref{ww11}, it holds
$$ \int_{\mathbb{S}^n } K \ov{u}  ^{ p-\tau}   \o  = \int_{\mathbb{S}^n } K (\a_0 \o)  ^{ p-\tau}   \o +  \sum_{j\geq 1} O\Big( \int \o  ^{p}  \wtilde{ \d }_j + \int \o  \wtilde{ \d }_j ^{p}  \Big) =
\a_0   ^{ p} \| \o \| ^2 + O\Big( \tau + \sum_{j\geq 1}  \frac{ 1 }{ \l_j ^{(n-2)/2} } \Big) . $$
Hence the proof  follows.
\end{pf}

\begin{pfn}{ \bf of Theorem \ref{th:blowupConv}} Let  $(\a_0,\a,\l,x,0)\in M_\tau$.  The proof goes along  the ideas introduced in \cite{BLR}.
We need to solve the system $((E_{\a_0}), (E_\a), (E_\l), (E_x), (E_v))$ introduced in Proposition \ref{G1}. \\
Proposition \ref{ovvmin} implies that $\ov{v}$ is a critical point of  $ I_{K,\tau}$ in the space $E_{\o, x, \l}^\perp$. Thus we get
$$ \n  I_{K,\tau}  \Big(\a_0 \o + \sum_{i=1}^{N} \a_i \wtilde{\d}_{ x_i,\l_i }+ \ov{v}\Big) _{\lfloor {E_{\o, x, \l}^\perp}} = 0 $$
and therefore (by the Lagrange multiplier theorem) there exist some constants $(A_0, A, B, C) \in \R \times \R^N \times \R^N \times (\R^n)^N$ such that
$$ \n  I_{K,\tau}  \Big(\a_0 \o + \sum_{i=1}^{N} \a_i \wtilde{\d}_{ x_i,\l_i }+ \ov{v}\Big)  = A_0 \o + \sum_{i=1}^N \Big( A_i \wtilde{\d}_{ x_i,\l_i } + B_i \l_i \frac{\partial \wtilde{\d}_{ x_i,\l_i }}{\partial \l_i } + C_i \cdot \frac{1}{\l_i} \frac{\partial \wtilde{\d}_{ x_i,\l_i }}{\partial x_i } \Big)$$
which implies that the equation $(E_v)$ in Proposition \ref{G1} is satisfied.

 Next, we estimate the numbers $A_0,A,B,C$ by taking the scalar product of $(E_v)$ with $\o$, $\wtilde{\d}_i$, $\l_i {\partial \wtilde{\d}_i}/{\partial \l_i}$ and  $\l_i^{-1}{\partial \wtilde{\d}_i}/{\partial x_i}$  respectively. Thus  we derive  a quasi-diagonal system in the variables $A_0$,  $A$, $B$ and $C_i$'s. The right  hand side is given by (using Propositions \ref{p:35w}, \ref{p:33w}, \ref{p:38w} and \ref{p:1conv}, the fact that $(\a_0, \a,\l,x,0)\in M_\tau$ and ${\partial \Psi_{\tau}}/{\partial v}= \n I_\tau (u)$ with $u:= \a_0 \o + \sum_{i=1}^{N}\a_i \wtilde{\d}_i  +\ov{v}$)
\begin{align*}
 & \langle  \frac{\partial \Psi_{\tau}}{\partial v}, \o \rangle = O (\tau \ln ^2\tau)  \quad ; \quad \quad    \langle  \frac{\partial \Psi_{\tau}}{\partial v}, \wtilde{\d}_i\rangle = O (\tau \ln ^2\tau) ; \\
  &    \langle \frac{\partial \Psi_{\e}}{\partial v}, 
 \l_i \frac{\partial \d_i}{\partial \l_i}\rangle = O (\tau )  \quad ; \quad \quad  \langle \frac{\partial \Psi_{\tau}}{\partial v},\frac{1}{\l_i}\frac{\partial \wtilde{\d}_i}{\partial x_i}\rangle =  O (\tau \ln ^2\tau) .
  \end{align*}
 Hence we deduce that
\be\label{F09} A_0=O ( \tau \ln^2\tau ) \, \,  ; \quad A_i=O ( \tau \ln^2\tau ) \, \,  ; \quad B_i=O (\tau ) \, \,  ;\quad  C_i=O ( \tau \ln^2\tau )  \quad \mbox{ for } i=1,\cdots,N.\ee
Furthermore, since $(\a_0, \a,\l,x,0) \in M_\tau$ then Proposition \ref{ovvmin} implies that $\| \ov{v}\| = O(\tau)$ and  therefore the system $((E_{\a_0}),(E_{\a_i}),(E_{\l_i}),(E_{x_i}))$ (introduced in Proposition \ref{G1}) is  equivalent to 

$$ (S): \qquad \begin{cases}
&  \langle \n  I_{K,\tau} (u),\o \rangle  =  0  \\
&  \langle \n  I_{K,\tau} (u),\wtilde{\d}_i \rangle  =  0  \quad \mbox{ for } i = 1,\cdots,N;\\
&  \langle \n  I_{K,\tau} (u),{\l_i}{\partial\wtilde{\d}_i}/{\partial \l_i}\rangle = O ( \tau^2 \ln^2\tau) \quad \mbox{ for } i = 1,\cdots,N; \\
& \langle \n  I_{K,\tau} (u),{\l_i}^{-1}{\partial\wtilde{\d}_i}/{\partial x_i}\rangle = O ( \tau^2 \ln^2\tau)  \quad \mbox{ for } i = 1,\cdots,N.
\end{cases} $$
\noindent
Now we introduce  the following change of variables:
\begin{align}
&  \b_0:= 1 - \a_0^{4/(n-2)} \, , \\
& \b_i:= 1 - \a_i^{4/(n-2)} K(x_i)\, , \quad i=1,\cdots,N;\label{b}\\
& \frac{1}{\l_i ^2} := - \frac{ c_4 K(y_i) }{ c_3  (\D K (y_i)} \tau (1+ \Lambda_i) \, , \quad i=1,\cdots,N;\label{l}\\
& x_i:= \frac{y_i+{\xi}_i}{ | y_i +{\xi}_i | } \quad \mbox{with } {\xi}_i \in T_{y_i}( \mathbb{S}^n) \, , \quad i=1,\cdots,N . \label{x}
\end{align}
Using this change of variables and Propositions \ref{p:35w}, \ref{p:33w}, \ref{p:38w} and \ref{p:1conv},  the previous system $(S)$ becomes
\be\label{S'1} (S') \, \begin{cases}
 \b_i = O ( \tau\ln\tau) \, ,  \quad \mbox{ for } i = 0,1,\cdots, N, \\
 \Lambda_i =   O ( \tau^{1/5} ) \, ,  \quad \mbox{ for } i = 1, \cdots,N, \\
D^2 K ( y_i) ( \xi_i ,.) = O (| \xi ^2 |  + \tau ^{{2}/{n}} \ln^2\tau+  \tau ^{{1}/{5} }) \, , \mbox{ for } i = 1, \cdots , N .
\end{cases} \ee
Since the critical points $y_i$'s are assumed to be non-degenerate, then, using the fixed point theorem, we derive the existence of $( \b_0^\tau, \b^\tau, \xi^\tau, \Lambda^\tau)$ such that the system $(S')$ is satisfied. Hence the existence of a critical point of $I_{K,\tau}$ (by Proposition \ref{G1}).\\
Note that, for $\o=0$, the unicity of such a solution is proved in \cite{MM}. The same argument holds for $\o\neq 0$ non-degenerate.\\
Lastly, for the Morse index, the same argument used in the proof of Theorem \ref{th:blowup} holds  which completes the proof of Theorem \ref{th:blowupConv}.
\end{pfn}

\section{Appendix}

In this appendix we collect various estimates needed through the paper.
 \begin{lem}\label{lowerL2}\cite{AB20b}
 Let $a\in { \mathbb{S}^n}$ and $\l > 0$ be large.\\
  $(i)$ Assume that $\tau \ln \l$ is small enough, then it holds
 \begin{align} \label{lower2}\d_{a,\l}^{-\tau}(x) = & c_0^{-\tau} \l^{-\tau (n-2)/2}\Big(1 + \frac{n-2}{2}\, \tau \ln (2+(\l^2-1)(1 -\cos d(a,x))) \Big) \\
 & + O\Big(\tau ^2 \ln (2+(\l^2-1)(1 -\cos d(a,x)))\Big) \quad \mbox{ for each } y \in {\mathbb{S}^n}. \nonumber \end{align}
 $(ii)$ For each $\g > 0$ and each $\b\in [0, n/(n-2))$, it holds
 $$ 0 < \int_{\mathbb{S}^n} \d_{a,\l}^{p+1 - \b } (x) \ln^\g \Big(2+(\l^2-1)(1 -\cos d(a,x))\Big)dx = O\Big( \frac{1}{\l^{\b (n-2)/2} }\Big).$$
 \end{lem}

Next, we state some elementary estimates needed in the paper. Their proofs follow from Taylor's expansion.

\begin{lem}
\begin{enumerate}
  \item[i)]
  Let $t_i > 0$ and $a,b \in \R$, there hold
\be \label{ti5}
 | (\sum t_i)^\g - \sum t_i ^\g | \leq c \begin{cases} \sum_{i\neq j} (t_it_j)^{\g/2} \quad & \mbox{ if }  0 < \g \leq 2 \\
 \sum_{i\neq j} t_i ^{\g-1} t_j  \quad & \mbox{ if }  \g > 2\end{cases} . \ee
\be \label{ti6}
 | | a+b|^\g - |a| ^\g - \g | a |^{\g-2} a b  | \leq c \begin{cases}  | b | ^\g + | a |^{\g-2} b^2 \quad & \mbox{ if }  \g > 2 \\
 | b | ^\g  \quad & \mbox{ if }  1 < \g \leq 2\end{cases} . \ee
  \item[ii)]
   For $ 1 < \g \leq 3$ and $t_1, \cdots, t_{N+1} > 0$, it holds
\be\label{est11} | (\sum t_i)^\g - \sum t_i^\g - \g t_1^{\g-1} (\sum_{j\neq 1} t_j) | t_1 \leq c \sum_{k\neq j } (t_k t_j)^{(\g + 1)/2} .\ee
  \item[iii)]
   For $ \g > 1 $, and $ x , y > 0 $, it holds
\be \label{est111}  ( x + y )^\g = x ^\g + O( x^{\g-1} y + y^\g ). \ee
  \item[iv)]
  For $ \g < 1 $, and $ x , y > 0 $, it holds
\be \label{est112}  ( x + y )^\g x = x ^{\g+1} + O( (xy)^{(\g+1)/2}  ). \ee
\end{enumerate}
\end{lem}
Finally, we prove a crucial estimate of the derivative of the infinite dimensional variable with respect to the parameter $ \a$.
\begin{lem}\label{Lemmalast}
The function $ v_\tau$,  defined in Proposition \ref{lambdaepsilonw}, satisfies 
$$
\Big| \Big|\frac{\partial v_\tau}{\partial \alpha}\Big|\Big| = o(1) . 
$$
\end{lem}
\begin{pf} Recall that we are in the case where $ u_\tau = \sum_{ i=1 }^N \a_i \wtilde{\d}_{a_i, \l_i} + \a \o + v_\tau $. \\
First, we claim that $ \partial v_\tau / \partial \a \in E_{\o, a, \l } ^\perp $ . \\
Indeed, since $ v_\tau  \in E_{\o, a, \l } ^\perp $, it follows that $ \langle \wtilde{\d}_{a_i, \l_i} , v_\tau \rangle = 0 $ for each $ \a $. Therefore, taking the derivative with respect to $ \a $, we deduce that $ \langle \wtilde{\d}_{a_i, \l_i} , \partial v_\tau / \partial \a  \rangle = 0 $. The other orthogonality constraints can be proved in the same way. Hence our claim follows. \\
Second, note that $ u _ \tau $ satisfies the problem $ ( \mathcal{N}_{K,\tau} ) $ and therefore the function $ v _ \tau $ satisfies the following PDE 
$$ {L}_{g_0} v_\tau := \Big( - \Delta  + \frac{n(n-2)}{4}  \Big) v_\tau =  K  u_\tau ^{p-\tau } - \sum_{i=1}^N \a_i \wtilde{\d }_{a_i, \l_i} ^p - \a  K \o^p $$
which implies that 
\be\label{val1}  \Big( - \Delta  + \frac{n(n-2)}{4}  \Big) \frac{\partial v_\tau }{\partial \a } =  (p-\tau) K  u_\tau ^{p-\tau-1 } \Big( \o + \frac{\partial v_\tau }{\partial \a }\Big) -   K \o^p . \ee
Following the proof of Lemma \ref{ovvmin} and using \eqref{val}, we write
\be \label{f1} 
 \frac{\partial v_\tau }{\partial \a } = \Big(  \frac{\partial v_\tau }{\partial \a } \Big)_- + \Big(  \frac{\partial v_\tau }{\partial \a } \Big)_0 + \Big(  \frac{\partial v_\tau }{\partial \a } \Big)_+ . \ee
Since $  {\partial v_\tau } / {\partial \a } \in E_{\o, a , \l }^\perp $, from \eqref{0q0}, we deduce that 
 \be \label{f2} \Big \| \Big(  \frac{\partial v_\tau }{\partial \a } \Big)_0 \Big \| = o  \Big(\Big \|  \frac{\partial v_\tau }{\partial \a }\Big \|  \Big)  . \ee
Concerning the other components, for example $  ( {\partial v_\tau } / {\partial \a } )_+ $, multiplying \eqref{val1} by $ ( {\partial v_\tau } / {\partial \a } )_+ $ and integrating over $ \mathbb{S}^n  $, we get 
\be\label{val2} \Big \| \Big( \frac{\partial v_\tau }{\partial \a } \Big)_+  \Big \| ^2 = (p-\tau) \int _{\mathbb{S}^n } K  u_\tau ^{p-\tau-1 }  \Big( \o + \frac{\partial v_\tau }{\partial \a } \Big) \Big( \frac{\partial v_\tau }{\partial \a } \Big)_+ - \int _{\mathbb{S}^n } K \o ^p \Big( \frac{\partial v_\tau }{\partial \a } \Big)_+ .
\ee
Note that, since $  ( {\partial v_\tau } / {\partial \a } )_+ \in N_+(\o) $, it follows that the  last integral in \eqref{val2} is equal to zero. For the other integral, we have 
\begin{align}
 \int _{\mathbb{S}^n } & K  u_\tau ^{p-\tau-1 } \o  \Big( \frac{\partial v_\tau }{\partial \a } \Big)_+ \nonumber \\
 &  =  \a^{p-\tau-1} \int _{\mathbb{S}^n } K  \o ^{p-\tau } \Big( \frac{\partial v_\tau }{\partial \a } \Big)_+ + \sum_{i=1}^N O \Big(   \int _{\mathbb{S}^n }   \Big\{\o ^{p -1 } \wtilde{\d}_{ a_i, \l_i}  + \wtilde{\d}_{ a_i, \l_i} ^{p - 1 } \o \Big\} \Big| \Big( \frac{\partial v_\tau }{\partial \a } \Big)_+  \Big|  \Big) \nonumber \\
 & = o \Big(  \Big \|   \Big( \frac{\partial v_\tau }{\partial \a }  \Big)_+ \Big \|  \Big)  \label{f3} \end{align}
 and 
\begin{align} \label{f4} 
(p-\tau)&  \int _{\mathbb{S}^n } K  u_\tau ^{p-\tau-1 }  \frac{\partial v_\tau }{\partial \a } \Big( \frac{\partial v_\tau }{\partial \a } \Big) _+ \nonumber \\
 & =  p \sum_{i=1}^N  \int _{\mathbb{S}^n }   \wtilde{ \d}_{a_i, \l_i }  ^{p - 1 }  \frac{\partial v_\tau }{\partial \a } \Big( \frac{\partial v_\tau }{\partial \a } \Big) _+  
 + p \int _{\mathbb{S}^n } K  \o ^{p -1 }  \frac{\partial v_\tau }{\partial \a } \Big( \frac{\partial v_\tau }{\partial \a } \Big) _+ + o \Big( \Big \|   \frac{\partial v_\tau }{\partial \a }  \Big\| \Big \|  \Big( \frac{\partial v_\tau }{\partial \a } \Big)_+ \Big\|  \Big).  
\end{align}
Since  $ H_0(\o) $, $ N_-(\o) $ and $ N_+ (\o) $ are orthogonal spaces for $ \langle . , \rangle $ and the bilinear form $ \int K \o^{p-1} .. $, we derive that 
\be \label{f5} 
\int _{\mathbb{S}^n } K  \o ^{p -1 }  \frac{\partial v_\tau }{\partial \a } \Big( \frac{\partial v_\tau }{\partial \a } \Big) _+ = \int _{\mathbb{S}^n } K  \o ^{p -1 }  \Big( \frac{\partial v_\tau }{\partial \a } \Big) _+ ^2 . \ee
Now, since $ H_0(\o) $ and $ N_-(\o) $  are finite dimensional spaces, we get 
\begin{align} 
\int _{\mathbb{S}^n }  & \wtilde{ \d}_{a_i, \l_i }  ^{p - 1 }  \frac{\partial v_\tau }{\partial \a } \Big( \frac{\partial v_\tau }{\partial \a } \Big) _+ \nonumber \\
& = \int _{\mathbb{S}^n }   \wtilde{ \d}_{a_i, \l_i }  ^{p - 1 }   \Big( \frac{\partial v_\tau }{\partial \a } \Big) _+ ^2 + O  \Big( \Big( \Big \|  \Big( \frac{\partial v_\tau }{\partial \a } \Big)_0  \Big\|_\infty  + \Big \|  \Big( \frac{\partial v_\tau }{\partial \a } \Big)_- \Big\|_\infty  \Big) \int _{\mathbb{S}^n }   \wtilde{ \d}_{a_i, \l_i }  ^{p - 1 } \Big|  \Big( \frac{\partial v_\tau }{\partial \a } \Big)_+ \Big| \Big)  \nonumber \\
 & = \int _{\mathbb{S}^n }   \wtilde{ \d}_{a_i, \l_i }  ^{p - 1 }   \Big( \frac{\partial v_\tau }{\partial \a } \Big) _+ ^2 + o  \Big( \Big \|  \frac{\partial v_\tau }{\partial \a }   \Big\|   \Big \|  \Big( \frac{\partial v_\tau }{\partial \a } \Big)_+ \Big\|   \Big) . \label{f6}
 \end{align}

However, since $ \o $  is bounded and $\langle  \partial v_\tau /\partial \a  , \o \rangle = 0  $, the first integral in \eqref{val2} satisfies 
\begin{align*}
 \int _{\mathbb{S}^n } K  u_\tau ^{p-\e-1 } \o \frac{\partial v_\tau }{\partial \a } & =  \a^{p-\e-1} \int _{\mathbb{S}^n } K  \o ^{p-\e } \frac{\partial v_\tau }{\partial \a } + O \Big(  \sum_{i=1}^N \int _{\mathbb{S}^n }   \Big\{\o ^{p-\e-1 } \wtilde{\d}_{ a_i, \l_i}  + \wtilde{\d}_{ a_i, \l_i} ^{p-\e-1 } \o \Big\} \Big|  \frac{\partial v_\tau }{\partial \a } \Big|  \Big) \\
 & = o \Big( \Big\|   \frac{\partial v_\tau }{\partial \a }  \Big\|  \Big).\end{align*}
Thus, using \eqref{f3}-\eqref{f6}, the equation \eqref{val2} becomes
$$ \Big \| \Big( \frac{\partial v_\tau }{\partial \a } \Big)_+  \Big \| ^2 = p \int _{\mathbb{S}^n } K  \o ^{ p - 1 }  \Big( \frac{\partial v_\tau }{\partial \a } \Big)_+^2 + p \sum_{i=1}^N \int _{\mathbb{S}^n }   \wtilde{\d}_{a_i, \l_i} ^{ p - 1 }  \Big( \frac{\partial v_\tau }{\partial \a } \Big)_+^2 + o  \Big( \Big \|  \frac{\partial v_\tau }{\partial \a }   \Big\|   \Big \|  \Big( \frac{\partial v_\tau }{\partial \a } \Big)_+ \Big\|  +   \Big \|  \Big( \frac{\partial v_\tau }{\partial \a } \Big)_+ \Big\| \Big) 
$$
which implies, by using Lemma \ref{Qwalpositive}, that 
$$ c \Big \|  \Big( \frac{\partial v_\tau }{\partial \a } \Big)_+ \Big\| ^2 \leq Q_{\o, a , \l } \Big(  \Big( \frac{\partial v_\tau }{\partial \a } \Big)_+ \Big) =  o  \Big( \Big \|  \frac{\partial v_\tau }{\partial \a }   \Big\|   \Big \|  \Big( \frac{\partial v_\tau }{\partial \a } \Big)_+ \Big\|  +   \Big \|  \Big( \frac{\partial v_\tau }{\partial \a } \Big)_+ \Big\| \Big) $$
and therefore
\be \label{f7}
\Big \|  \Big( \frac{\partial v_\tau }{\partial \a } \Big)_+ \Big\| = o(1) + o  \Big( \Big \|  \frac{\partial v_\tau }{\partial \a }   \Big\|  \Big) . \ee
In the same way, multiplying \eqref{val1} by $ ( {\partial v_\tau } / {\partial \a } )_- $ and integrating over $ \mathbb{S}^n  $, we get 
\be \label{f8}
\Big \|  \Big( \frac{\partial v_\tau }{\partial \a } \Big)_- \Big\| = o(1) + o  \Big( \Big \|  \frac{\partial v_\tau }{\partial \a }   \Big\|  \Big) . \ee
Combining \eqref{f2}, \eqref{f7} and \eqref{f8}, we get that 
$$ \Big \|   \frac{\partial v_\tau }{\partial \a }  \Big\| = o(1) $$
which completes the proof of the lemma.
\end{pf}

\end{document}